\documentclass[11pt]{article}
\usepackage{graphicx, xypic, color, url}
\usepackage{amssymb}
\usepackage{epstopdf}
\usepackage{hyperref}
\usepackage{enumerate}

\input epsf.sty
\epsfverbosetrue 

\usepackage{latexsym}
\usepackage{amsfonts}
\usepackage{amssymb}

\newcommand{\R}{\mathbb{R}}

\newcommand{\Z}{\mathbb{Z}}
\newcommand{\C}{\mathbb{C}}
\newcommand{\IP}{\mathbb{P}}
\newcommand{\IH}{\mathbb{H}}

\newcommand{\N}{\mathbb{N}}
\newcommand{\D}{\mathbb{D}}

\newcommand{\rs}{\mbox{$\widehat{\C}$}}

\def\SSS{{\mathcal S}}
\def\TTT{{\mathcal T}}
\def\AAA{{\mathcal A}}

\def\FFF{{\mathcal F}}

\def\HHH{{\mathcal H}}

\def\MMM{{\mathcal M}}
\def\QQQ{{\mathcal Q}}
\def\OOO{{\mathcal O}}
\def\MMM{{\mathcal M}}
\def\WWW{{\mathcal W}}

\def\SSS{{\mathcal S}}

\def\ZZZ{{\mathcal Z}}

\newtheorem{thm}{Theorem}[section]
\newtheorem{defn}[thm]{Definition}

\newtheorem{prop}[thm]{Proposition}

\newtheorem{lemma}[thm]{Lemma}
\newtheorem{cor}[thm]{Corollary}

\newcommand{\qed}{\nopagebreak \begin{flushright}
        \rule{2mm}{2.5mm} \end{flushright}}

\newcommand{\implies}{\Rightarrow}
\newcommand{\bdry}{\partial}                     
\newcommand{\id}{\mbox{\rm id}}                  
\newcommand{\cl}{\overline}                      

\newcommand{\Aut}{\mbox{\rm Aut}}                        
\newcommand{\Mod}{\mbox{\rm Mod}}    
\newcommand{\Teich}{\mbox{\rm Teich}}
\newcommand{\intersect}{\cap}                    
\newcommand{\union}{\cup}                        

\newcommand{\mtwo}[4]                            
{\mbox{$\left(\begin{array}{cc}                  
#1 & #2 \\
#3 & #4 
\end{array}
\right)$}}

\newcommand{\dettwo}[4]                          
{\mbox{$\left|\begin{array}{cc}                  
#1 & #2 \\
#3 & #4 
\end{array}
\right|$}}

\newcommand{\pf}{\noindent {\bf Proof: }}

\newcommand{\be}{\begin{enumerate}}
\newcommand{\eb}{\end{enumerate}}
\newcommand{\bi}{\begin{itemize}}
\newcommand{\ib}{\end{itemize}}
\newcommand{\bl}{\begin{list}}
\newcommand{\lb}{\end{list}}

\newcommand{\gap}{\vspace{5pt}}                 

\newcommand{\Stab}{\mbox{\rm Stab}}

\newcommand{\Fix}{\mbox{\rm Fix}}

\newcommand{\Homeo}{\mbox{Homeo}}

\newcommand{\genby}[1]{\mbox{$\langle #1 \rangle$}}

\def\Pf{{P_f}}

\def\P{{{\mathbb P}}}

\def\Teich{\mathrm{Teich}}

\def\Tw{\mathrm{Tw}}
\def\PMod{\mathrm{PMod}}
\def\supp{\mathrm{supp}}

\newcommand{\ns}[1]{\textcolor{red}{ #1}}

\definecolor{orange}{rgb}{1,0.5,0}

\newcommand{\dom}{\mbox{dom}}
\newcommand{\pullback}{\stackrel{f}{\longleftarrow}}

\def\Pf{{P_f}}

\def\C{{\mathbb C}}
\def\D{{\mathbb D}}

\def\N{{\mathbb N}}
\def\0{{\mathbb 0}}
\def\P{{{\mathbb P}}}

\def\R{{\mathbb R}}

\def\Z{{\mathbb Z}}

\def \epsilon{\varepsilon}

\def\Teich{\mathrm{Teich}}

\def\id{{\rm id}}
\def\qed{{\hfill{$\square$}}}

\def\pmcg{\mathrm{PMCG}(\P^1, {\Pf})}
\newcommand{\Config}{\mbox{\rm Config}}

\begin{document}

\title{Pullback invariants of Thurston maps}

\author{Sarah Koch, Kevin M. Pilgrim, and Nikita Selinger}
\maketitle

\begin{abstract}
Associated to a Thurston map $f: S^2 \to S^2$ with postcritical set $P$ are several different invariants obtained via pullback: a relation $\SSS_{P} \pullback \SSS_{P}$ on the set $\SSS_{P}$ of free homotopy classes of curves in $S^2\setminus P$, a linear operator $\lambda_f: \R[\SSS_{P}]\to\R[\SSS_{P}]$ on the free $\R$-module generated by $\SSS_{P}$, a virtual endomorphism $\phi_f: \PMod(S^2, P) \dashrightarrow \PMod(S^2, P)$ on the pure mapping class group, an analytic self-map  $\sigma_f: \TTT(S^2, P) \to \TTT(S^2, P)$ of an associated Teichm\"uller space, and an analytic self-correspondence $X\circ Y^{-1}: \MMM(S^2, P) \rightrightarrows \MMM(S^2, P)$ of an associated moduli space. Viewing these associated maps as invariants of $f$, we investigate relationships between their properties. 
\end{abstract}

\newpage

\tableofcontents 

\section{Introduction}

Thurston maps are orientation-preserving branched covers $f: S^2 \to S^2$ from the oriented topological two-sphere to itself that satisfy certain properties. They were introduced by Thurston as combinatorial invariants associated to postcritically finite rational functions $R: \rs \to \rs$, regarded as dynamical systems on the Riemann sphere.  A fundamental theorem of complex dynamics is Thurston's characterization and rigidity theorem  \cite{DH1}, which (i) characterizes which Thurston maps $f$ arise from rational functions, and (ii) says that apart from a well-understood family of ubiquitous exceptions, the M\"obius conjugacy class of $R$ is determined by the combinatorial class of its associated Thurston map $f$.  The proof of this theorem transforms the question of determining whether $f$ arises from a rational function $R$ to the question of whether an associated {\em pullback map} 
$\sigma_f: \TTT_P \to \TTT_P$, an analytic self-map of a certain Teichm\"uller space, has a fixed point. 

As combinatorial (as opposed to analytic or algebraic) objects, Thurston maps, in principle, should be easier to classify than postcritically finite rational maps. For many years, the lack of suitable invariants for general Thurston maps frustrated attempts toward this goal.  In a ground-breaking paper \cite{bartholdi:nekrashevych:twisted}, Bartholdi and Nekrasheyvh introduced several new tools.  One innovation was to develop new algebraic tools from the theory of self-similar groups to study the dynamics of a given Thurston map $f:S^2\to S^2$.  Another was to exploit the existence of an associated analytic correspondence on the moduli space $\MMM_P$ covered by the graph of $\sigma_f$.   These new tools have led to much better invariants for Thurston maps and to a better understanding of the pullback map $\sigma_f$, see \cite{nekrashevych:combinatorics}, \cite{cfpp:fsr-cve}, \cite{koch:thesis}, \cite{koch:criticallyfinite}, \cite{kmp:tw}, \cite{bekp}, \cite{kmp:ph:cxcii},  \cite{kmp:ph:expanding}, \cite{kelsey:schemes}, \cite{nekrashevych:expanding}, \cite{bartholdi:nekrashevych:imgquad1}, \cite{lodge:thesis}. 

In this work, we deepen the investigations of the relationship between dynamical, algebraic, and analytic invariants associated to Thurston maps.   
\gap

\noindent{\bf Fundamental invariants.} 
Let $S^2$ denote a topological two-sphere, equipped with an orientation. Fix an identification $S^2=\IP^1=\rs$; we use $S^2$ for topological objects, $\IP^1$ for holomorphic objects, and $\rs$ when formulas are required. 
Let $P \subset S^2$ be a finite set with $\#P\geq 3$ (in the case $\#P=3$, the groups and spaces are trivial, so it is helpful to imagine at first that $\# P \geq 4$). 

The following objects are basic to our study. 
\be
\item $\SSS_P$, the set of free homotopy classes of simple, unoriented, essential, nonperipheral, closed curves in $S^2\setminus P$; the symbol $o$  denotes the union of  the sets of free homotopy classes of unoriented inessential and peripheral curves. A curve representing an element of $\SSS_P$ we call {\em nontrivial}.  A multicurve $\Gamma$ on $S^2\setminus P$ is a possibly empty subset of $\SSS_P$ represented by nontrivial, pairwise disjoint, pairwise  nonhomotopic curves. Let $\MMM\SSS_P$ denote the set of possibly empty multicurves on $S^2\setminus P$.
\item $\R[\SSS_P]$, the free $\R$-module generated by $\SSS_P$; this arises naturally in the statement of Thurston's characterization theorem.  
\item $G_P:=\PMod(S^2, P)$, the pure mapping class group (that is, orientation-preserving homeomorphisms $g: (S^2, P) \to (S^2, P)$ fixing $P$ pointwise, up to isotopy through homeomorphisms fixing $P$ pointwise). Each nontrivial element of $G_P$ has infinite order.  The group $G_P$ contains a distinguished subset $\Tw$ whose elements are multitwists;  that is 
\[
\Tw:=\bigcup_{\Gamma\text{ a multicurve on }S^2\setminus P} \Tw(\Gamma)
\]
where $\Tw(\Gamma)$ is the subgroup of $G_P$ generated by Dehn twists around components of $\Gamma$. A multitwist is {\em positive} if it is a composition of positive powers of right Dehn twists.  Let 
\[
\Tw^+:=\bigcup_{\Gamma\text{ a multicurve on }S^2\setminus P} \Tw^+(\Gamma)
\]
where $\Tw^+(\Gamma)$ is the subgroup of $G_P$ generated by positive Dehn twists around components of $\Gamma$. The {\em support} of a multitwist is the smallest multicurve\ns{,} about which the twists comprising it occur; this is is well-defined. 

\item $\TTT_P:=\Teich(S^2, P)$, the Teichm\"uller space of $(S^2, P)$, as in \cite{DH1}.  It comes equipped with two natural metrics, the Teichm\"uller Finsler metric $d_{\TTT_P}$ and the Weil-Petersson (WP) inner product metric $d^{WP}_{\TTT_P}$.  The metric space $(\TTT_P, d_{\TTT_P})$ is complete, and any pair of points are joined by a unique geodesic.  In contrast, $(\TTT_P, d^{WP}_{\TTT_P})$ is incomplete. The WP-completion $\cl{\TTT}_P$ has a rich geometric structure \cite{masur:duke:1976} (see also \cite{hubbard:koch:dm});  it is a stratified space where each  stratum $\TTT_{P}^\Gamma$  corresponds to a multicurve $\Gamma$ and consists of noded Riemann surfaces whose nodes correspond to pinching precisely those curves comprising $\Gamma$ to points.  We denote by $\bdry \TTT_P$ the WP boundary of $\TTT_P$. 

\item $\MMM_P:=\MMM(S^2, P)$ is the corresponding moduli space.   It is a complex manifold, isomorphic to a complex affine  hyperplane complement. The natural projection $\pi:\TTT_P \to \MMM_P$ is a universal cover, with deck group $G_P$. 
\eb 

Given a basepoint $\tau_\circledast \in \TTT_P$, there is a natural identification of $G_P$ with $\pi_1(\MMM_P, m_\circledast)$ where $m_\circledast:=\pi(\tau_\circledast)$; the isomorphism $\pi_1(\MMM_P, m_\circledast) \to G_P$ proceeds via isotopy extension, while the isomorphism $G_P \to \pi_1(\MMM_P, m_\circledast)$ is induced by composing the evaluation map at those points marked by $P$ with an isotopy to the identity through maps fixing three points of $\MMM_P$; see \cite{lodge:thesis} for details. 

 A subgroup $L<G_P$ is {\em purely parabolic} if $L \subset \Tw$. In terms of the identification $G_P \leftrightarrow \pi_1(\MMM_P, m_\circledast)$, $L$ is purely parabolic if and only if the displacements of each of its elements satisfies $\inf_{\tau \in \TTT_P}d_{\TTT_P}(\tau, g.\tau)=0$ for all $g\in L$. A purely parabolic subgroup is {\em complete}  if for each multitwist $M \in L$ and each $\gamma \in \SSS_P$, we have $\gamma \in \supp(M) \implies T_\gamma \in L$ where $T_\gamma$ is the right Dehn twist about $\gamma$. A complete parabolic  subgroup is necessarily of the form $\Tw(\Gamma)$, the set of all multitwists about elements of $\Gamma$, where $\Gamma$ is some nonempty multicurve. 

\gap

Let $f:S^2\to S^2$ be an orientation-preserving branched covering map with critical set $\Omega_f$. Its {\em postcritical set} is 
\[ P_f:= \bigcup_{n>0}f^{\circ n}(\Omega_f); \] 
we assume throughout this work that $\#P_f < \infty$.  Suppose $P \subset S^2$ is finite, $P_f \subset P$, and $f(P) \subset P$; we then call the map of pairs $f: (S^2, P) \to (S^2, P)$ a {\em Thurston map}.  Throughout this work, we refer to a Thurston map by the symbol $f$, and often suppress mention of the non-canonical subset $P$. 

The following objects are associated to a Thurston map $f:(S^2, P) \to (S^2, P)$ via pullback:
\be
\item a relation $\beta_f: \SSS_P \union \{o\} \to \SSS_P \union \{o\}$, 
\item a non-negative linear operator $\lambda_f: \R[\SSS_P] \to \R[\SSS_P]$, 
\item a virtual endomorphism $\phi_f: G_P \dashrightarrow G_P$, 
\item an analytic map $\sigma_f: \TTT_P\to\TTT_P$; this extends continuously to $\sigma_f: \cl{\TTT}_P \to \cl{\TTT}_P$ by \cite[\S 4]{selinger:boundary}, 
\item an analytic correspondence $X\circ Y^{-1}: \MMM_P \rightrightarrows \MMM_P$; the double-arrow notation reflects our view that this is a one-to-many ``function''.
\eb
Precise definitions and a summary of basic properties of these associated maps will be given in the next two sections.  
\gap 

\noindent{\bf Main results.}  The goal of this work is to examine how properties of the maps in $(1)$ through $(5)$ are related.  That such relationships should exist is expected since there are fundamental identifications between various objects associated to the domains of these maps; see \S 2.  Our main results include the following. 

\bi
\item We characterize when the map on Teichm\"uller space $\sigma_f: \TTT_P \to \TTT_P$ on Teichm\"uller space sends the Weil-Petersson boundary to itself  (Theorem \ref{thm:bdry}).  
\item We characterize when the map on Teichm\"uller space $\sigma_f: \TTT_P \to \TTT_P$ is constant (Theorem \ref{thm:constant}); the proof we give corrects an error in \cite{bekp}. 
\item We show that the virtual endomorphism $\phi_f: G_P \dashrightarrow G_P$ is essentially never injective (Theorem \ref{thm:noninjective}). 
\item It is natural to investigate whether a Thurston map $f$ induces an action on projectivized measured foliations, i.e. on the Thurston boundary of $\TTT_P$.  Using Theorem \ref{thm:noninjective}, however, we show (Theorem \ref{thm:dense}) that each nonempty fiber of the pullback map $\sigma_f$ accumulates on the entire Thurston boundary.  This substantially strengthens the conclusion of  \cite[Theorem~9.4]{selinger:boundary}.
\item We give sufficient analytic criteria (Theorem \ref{thm:eventually_finite}) for the existence of a finite global attractor for the pullback relation $\SSS_P \pullback \SSS_P$ on curves---equivalently, for the action of the pullback map $\sigma_f: \cl{\TTT}_P \to \cl{\TTT}_P$ on strata. 
\item We prove an orbit lifting result (Proposition \ref{prop:shadowing}) which asserts that under the hypothesis that the virtual endomorphism $\phi_f$ is surjective, finite orbit segments of the pullback correspondence $X\circ Y^{-1}$ on moduli space can be lifted to finite orbit segments of $\sigma_f$. 
\item In the case when $\#P=4$, we relate fixed-points of the various associated maps (Theorem \ref{thm:tfae}).  If in addition the inverse of the pullback correspondence is actually a function (i.e. the map $X$ is injective),  sharper statements are possible (Theorem \ref{thm:dyncons}). Our results generalize, clarify, and put into context the algebraic, analytic, and dynamical findings in the analysis of twists of $z\mapsto z^2+i$ given in \cite[\S 6]{bartholdi:nekrashevych:twisted}.  Using the shadowing result from Proposition \ref{prop:shadowing}, these results also demonstrate that for certain unobstructed Thurston maps, one can build finite orbits of the pullback map whose underlying surfaces behave in prescribed ways.  For example, one can arrange so that the length of the systole can become shorter and shorter for a while before stabilizing.  
\ib

While motivated by the attempt at a combinatorial classification of dynamical systems, many of the results we obtain are actually more naturally phrased for a nondynamical branched covering map $f: (S^2, A) \to (S^2, B)$.  Where possible, we first phrase and prove more general results, (Proposition \ref{virtual}, Theorem \ref{thm:bdry}, Theorem \ref{thm:constant}, Proposition \ref{BvsA}, Theorem \ref{thm:noninjective}, and Theorem \ref{thm:dense}); the aforementioned theorems then become corollaries of these more general results. 
\gap

\noindent{\bf Conventions and notation.} $S^2$ denotes the topological two-sphere, equip\-ped with the orientation induced from the identification with $\rs$ via the usual stereographic projection. All branched covers and homeomorphisms are orientation-preserving. Throughout, $f$ denotes a branched covering of $S^2\to S^2$ of degree $d \geq 2$. The symbols $A, B, P$ denote finite subsets of $S^2$, which contain at least three points. 
\gap

\noindent{\bf Acknowledgements.} We thank Indiana University for supporting the visits of S. Koch and N. Selinger. S. Koch was supported by an NSF postdoctoral fellowship and K. Pilgrim by a Simons collaboration grant.  We thank A. Edmonds, D. Margalit, and V. Turaev for useful conversations.  
\gap

\section{Fundamental identifications}

\subsection{Preliminaries}

\noindent{\bf Teichm\"uller and moduli spaces.}   The Teichm\"uller space $\TTT_P=\TTT(S^2,P)$ is the space of equivalence classes of orientation-preserving homeomorphisms $\varphi:S^2\to \mathbb{P}^1$, whereby $\varphi_1\sim \varphi_2$ if there is a M\"obius transformation $\mu:\mathbb{P}^1\to\mathbb{P}^1$ so that 
\begin{itemize}
\item $\varphi_1=\mu\circ \varphi_2$ on the set $P$, and 
\item $\varphi_1$ is isotopic to $\mu\circ \varphi_2$ relative to the set $P$. 
\end{itemize}
The moduli space $\MMM_P=\MMM(S^2,P)$ is the set of all injective maps $\varphi:P\rightarrowtail\mathbb{P}^1$ modulo postcomposition by M\"obius transformations. The Teichm\"uller space and the moduli space are complex manifolds of dimension $\#P-3$, and the map $\pi_P:\TTT_P\to \MMM_P$ given by $\pi_P:[\varphi]\mapsto [\varphi|_P]$, is a holomorphic universal covering map.  

Note that since we have identified $S^2=\mathbb{P}^1$, both $\MMM_P$ and $\TTT_P$ have natural basepoints represented by the classes of the inclusion and identity maps, respectively. 
\gap

\noindent{\bf Teichm\"uller metric.} Equipped with the Teichm\"uller Finsler metric $d_{\TTT_P}$, the space $\TTT_P$ becomes a complete uniquely geodesic metric space, homeomorphic to the open ball $B^{\#P-3}$; it is not, however, nonpositively curved. 
\gap

\noindent{\bf WP metric.} In contrast, when equipped with the WP metric $d^{WP}_{\TTT_P}$, the space $\TTT_P$ is negatively curved but incomplete.  The completion $\cl{\TTT}_P$ is a stratified space whose strata  are indexed by (possibly empty) multicurves.  Each stratum is homeomorphic to the product of the Teichm\"uller spaces of the  components of the noded surface obtained by  collapsing each curve of $\Gamma$ to a point.   Indeed, this completion inherits an analytic structure and coincides with the augmented Teichm\"uller space parametrizing noded Riemann surfaces marked by $P$ (see \cite{hubbard:koch:dm}).  It is noncompact,  coarsely negatively curved, and quasi-isometric to the pants complex; for an extensive overview, see \cite{wolpert:handbook}. 
\gap 

\noindent {\textbf{Thurston compactification.}} For $\gamma \in \SSS_P$ and $\tau \in \TTT_P$ let $\ell_\gamma(\tau)$ denote the length of the unique geodesic in the hyperbolic surface associated to $\tau$. The map $\tau \mapsto (\ell_\gamma(\tau))_{\gamma \in \SSS_P} $ defines an embedding $\TTT_P \to \R_{\geq 0}^{\SSS_P}$ which projects to an embedding $\TTT_P \to \P\R_{\geq 0}^{\SSS_P}$.  Thurston showed that the closure of the image is homeomorphic to the closed ball $\cl{B}^{\#P-3}$, and that the boundary points may be identified with {\em projective measured foliations} on $(S^2, P)$. 
A comprehensive reference is the book by Ivanov \cite{ivanov:book:subgroups}.
\gap  

Fix a basepoint $\tau_\circledast \in \TTT_P$; this gives rise to a  basepoint $m_\circledast:=\pi(\tau_\circledast)\in \MMM_P$; recall that we then have an identification $G_P=\mathrm{PMod}(S^2,P)\leftrightarrow \pi_1(\MMM_P, m_\circledast)$. The following folklore theorem is well-known. 

\begin{thm}[Fivefold way]
\label{thm:fivefold_way}
There are natural bijections between the following sets of objects:
\be
\item multicurves $\Gamma$ on $S^2\setminus P$
\item  ``purely atomic'' measured foliations $\FFF(\Gamma):=\sum_{\gamma \in \Gamma} \nu_\gamma$, where $\nu_\gamma$ is the delta-mass at $\gamma$, which counts the number of intersections of a curve with $\gamma$. 
\item complete purely parabolic subgroups $L$ of $G$;  
\item strata $\TTT_P^\Gamma \subset \cl{\TTT}_P$;
\item certain subgroups of loops in moduli space (thought of as generated by certain pure braids) and corresponding via ``pushing'' to complete purely parabolic subgroups $L$ of $\pi_1(\MMM_P, m_\circledast)$.
\eb
\end{thm}
The subgroups arising in (5) will be described shortly. 
The bijections are given as follows: 

$(1) \rightarrow (2)$ Take the zero measured foliation if $\Gamma$ is empty; otherwise: take a foliation of a regular neighborhood of $\Gamma$ with the width (transverse measure) of each neighborhood equal to $1$. 

$(2) \rightarrow (1)$ If the foliation is the zero foliation, take $\Gamma = o$, the empty multicurve; otherwise, take one core curve from each cylinder in the normal form of the foliation.

$(1) \rightarrow (3)$ Take the subgroup $L:=\Tw(\Gamma)$ generated by Dehn twists about elements of $\Gamma$. 

$(3) \rightarrow (1)$  Take the union $\Gamma$ of the core curves of the representing twists; this is well-defined. 

$(1) \rightarrow  (4)$ We take those marked noded spheres in which precisely all the curves comprising $\Gamma$ correspond to nodes, and nothing else. 

$(3) \rightarrow (4)$ The stratum $\TTT_P^\Gamma$ is the unique stratum in $\Fix(L)$ of maximum dimension. 

$(4) \rightarrow (3)$ Given the stratum $\TTT_P^\Gamma$, take the pointwise stabilizer $L:=\Tw(\Gamma)=\Stab_{G_P}(\TTT_P^\Gamma)$.

$(3)\rightarrow (5)$ We now describe this correspondence. 
\gap

Let $P:= \{p_1, \ldots, p_n\} \subset S^2$, $n:=\#P \geq 4$. Recall that the {\em configuration space} $\Config(S^2, P)$ is the space of injections $P \rightarrowtail S^2$; the inclusion $P \subset S^2$ gives a natural basepoint. Let $\pi: \Config(S^2, P) \to \MMM_P$ be the natural projection and $\pi_*$ the induced map on fundamental groups. 
It is well-known that $\pi_1(\Config(S^2, P)) \simeq PB_{n-1}$, the pure Artin braid group on $n-1$ strands. 
Let $\Gamma$ be a nonempty multicurve represented by curves $\gamma \subset S^2\setminus P$. 
Below, we describe a procedure for producing, for each $\gamma \in \Gamma$, a loop $\ell_\gamma$ in $\Config(S^2, P)$ based at $P$ which projects to the right Dehn twist $\Tw(\gamma)$. This procedure will have the additional  property that as elements of $\pi_1(\Config(S^2, P))$, given $\gamma_1, \gamma_2 \in \Gamma$, the elements represented by $\ell_{\gamma_1}, \ell_{\gamma_2}$ will commute. 

The idea is based on simultaneous ``pushing'' (or Birman spin) of points; see Figure 1.   Suppose $z_1, \ldots, z_m \in \D$ are nonzero complex numbers, pick $r$ with $\max_i \{|z_i|\} < r < 1$, set $z_0=0$, and consider the motion $t \mapsto z_i(t):=\exp(2\pi i t)z_i, t \in [0,1], i=0, \ldots, m$. This motion extends to an isotopy of $\cl{\D}$  fixing the origin and the boundary.  Taking the result of this extension when $t=1$, the resulting ``multi-spin'' homeomorphism of the plane is homotopic, through homeomorphisms fixing $z_0, z_1, \ldots, z_m$, to the right-hand Dehn twist about the circle $\gamma:=\{|z|=r\}$ surrounding the $z_i$'s. To see this, note that the extension of each individual motion $z_i(t)$ yields a ``spin'' which is the composition of a  left Dehn twist about the left-hand component of a regular neighborhood of the circle traced by $z_i(t)$ with the right Dehn twist about the corresponding right-hand component. The left twist resulting from the motion of $z_1(t)$ is trivial, since this left boundary component is peripheral, and for $i=1, \ldots, m-1$ the right twist from the motion of $z_i$ cancels the left twist from the motion of $z_{i+1}$. 
\begin{figure}
\begin{center}
\includegraphics[width=3in]{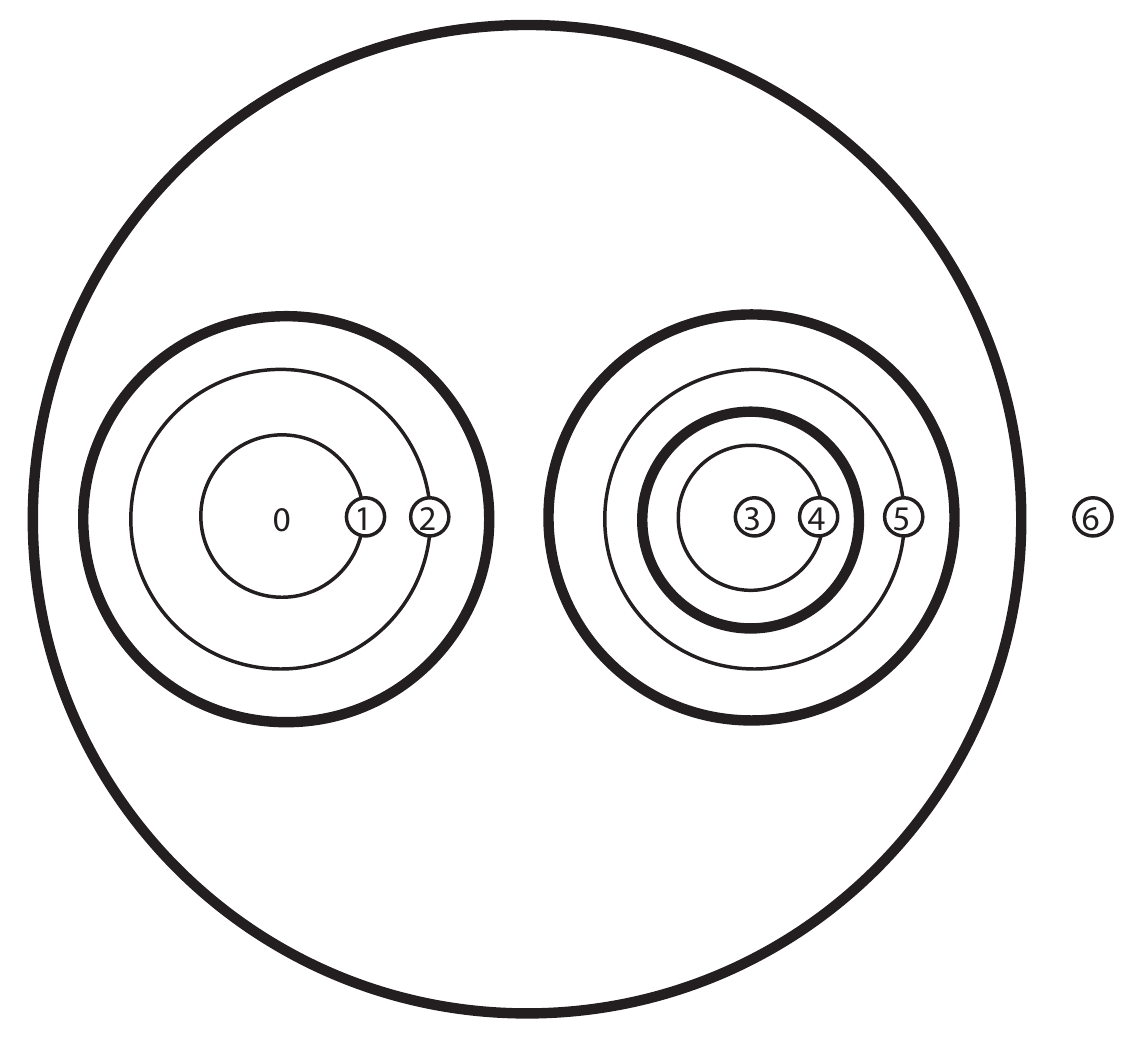}
\end{center}
\caption{Simultanteously pushing the points labelled 4 and 5 around the point labelled 3 yields the right Dehn twist about the bold curve on the outer side of the thin curve passing through point 5. }
\end{figure}

Now suppose $\Gamma$ is a nonempty multicurve. Choose an element $\delta \in \Gamma$  among the possibly several components of $\Gamma$ such that $\delta$ does not separate a pair of elements of $\Gamma$. Let $V$ be a component of $S^2\setminus \delta$ whose closure contains $\Gamma$ (the component $V$ is unique if $\#\Gamma > 1$).  Pick $\gamma \in \Gamma$. Then $\gamma$ bounds a distinguished Jordan domain $D_\gamma \subset V$.  Let  $\{p_{j_0}, p_{j_1}, \ldots, p_{j_m}\} = D_\gamma \intersect P$ where $1 \leq j_0 < j_1 < \ldots < j_m \leq n$, so that $j_0$ is the smallest index of an element of $P$ in $D_\gamma$. Up to postcomposition with rotations about the origin, there is a unique Riemann map $\phi: (D_\gamma, p_{j_0}) \to (\D, 0)$. Set $z_i =\phi(p_{j_i}), i=0, \ldots, m$, and transport the motion of the $z_i$ constructed in the previous paragraph to a motion of $\{p_{i_0}, p_{i_1}, \ldots, p_{i_m}\}$ in $S^2$; it is supported in the interior of $D_\gamma$.  This motion gives a loop $\ell_\gamma$ in the space  $\Config(S^2, P)$. 

By construction, $\ell_\gamma$ projects to the right Dehn twist about $\gamma$ in the pure mapping class group.  
Suppose $\gamma_1, \gamma_2$ are distinct elements of $\Gamma$. If $D_{\gamma_1} \intersect D_{\gamma_2}=\emptyset$ then clearly the elements represented by $\ell_{\gamma_1}$ and $\ell_{\gamma_2}$ commute. Otherwise, we may assume $D_{\gamma_1} \subset D_{\gamma_2}$. By construction, the loop $\ell_{\gamma_2}$ represents the central element in the braid group on $\#P \intersect D_{\gamma_2}$ strands, so it commutes with the element represented by $\ell_{\gamma_1}$. 

\gap

\section{Pullback invariants}

Many of the objects we are concerned with arise nondynamically. We first define them in the nondynamical setting (\S 3.1), and then consider the same objects in the dynamical case (\S 3.2). 

\subsection{No dynamics} 

\noindent{\bf Admissible covers.} Suppose $A, B \subset S^2$ are finite sets, each containing at least three points. A branched covering $f: (S^2, A) \to (S^2, B)$ is {\em admissible} if (i) $B \supseteq V_f$, the set of branch values of $f$, and (ii) $f(A) \subseteq B$; we do not require $A=f^{-1}(B)$. 

For the remainder of this subsection, we fix an admissible cover $f:(S^2,A)\to (S^2,B)$.  To $(f, A, B)$ we associate the following objects.  Though they depend upon all three elements of the triple, we indicate this dependency by referring only to $f$, for brevity. 
\gap

\noindent{\bf Pullback relation on curves.} The {\em pullback relation} 
\[
\SSS_B\cup \{o\}\pullback \SSS_A\cup\{o\}
\]
is defined by setting $o \pullback o$ and 
\[ \gamma_1 \pullback \gamma_2\]
if and only if  $\gamma_2$ is homotopic in $S^2\setminus A$ to a connected component of the preimage of $\gamma_1\subset S^2\setminus B$ under $f$. Thus $\gamma \pullback o$ if and only if some preimage of $\gamma$ is inessential or peripheral in $S^2\setminus A$.   The pullback relation induces a {\em pullback function} $f^{-1}: \MMM\SSS_B  \to \MMM\SSS_A $, by sending the empty multicurve on $(S^2, A)$ to the empty multicurve on $(S^2, B)$  and $\Gamma \mapsto f^{-1}(\Gamma):=\{ \tilde{\gamma} \in \SSS_A: \exists \gamma \in \Gamma, \gamma \pullback \tilde{\gamma}\}$; note that $f^{-1}(\Gamma)$ might be empty.  

\gap 

\noindent{\bf The associated linear transformation.} There is a linear transformation 
\[
\lambda_f:\R[\SSS_B]\to \R[\SSS_A]
\] 
defined on basis vectors by 
\[ 
\lambda_f(\gamma) = \sum_{\gamma \pullback \gamma_i} d_i \gamma_i\quad\text{where}\quad d_i:= \sum_{f^{-1}(\gamma)\supset \delta \simeq \gamma_i}\frac{1}{\deg(f:\delta \to \gamma)}.
\]
\noindent{\bf The pullback map.} There is a map $\sigma_f:\TTT_B\to\TTT_A$ associated to an admissible cover $f:(S^2,A)\to (S^2,B)$. Let $\varphi:(S^2,B)\to(\P^1,\varphi(B))$ be an orientation-preserving homeomorphism. Then there is an orientation-preserving homeomorphism $\psi:(S^2,A)\to (\P^1,\psi(A))$, and a rational map $F:(\P^1,\psi(A))\to(\P^1,\varphi(B))$ so that the following diagram commutes;
\[
\xymatrix{
(S^2,A)\ar[r]^\psi\ar[d]^f & (\P^1,\psi(A))\ar[d]^F\\
(A^2,B)\ar[r]^\varphi & (\P^1,\varphi(B))}
\]
we set $\sigma_f:[\varphi]\mapsto [\psi]$. It is well-known that this map is well-defined, holomorphic,  and distance-nonincreasing for the corresponding Teichm\"uller metrics on the domain and range \cite{DH1}. 

\bigskip

\noindent{\bf The virtual homomorphism.}   
\begin{prop}\label{virtual}
Let $f: (S^2, A) \to (S^2, B)$ be an admissible branched covering. Then the subset 
\[ H_f:=\{[h] : \exists \ \tilde{h}, \ h \circ f = f \circ \tilde{h}, \ \mbox{and} \ \tilde{h}|_A=\id_A\}\]
is a finite-index subgroup of $G_B$, and the function 
\[ \phi_f: H_f \to G_A, \;\; [h] \mapsto [\tilde{h}]\]
induced  by lifting representatives is well-defined and a homomorphism.
\end{prop}
\pf The possibility that the unramified covering induced by $f$ might admit deck transformations complicates the proof. 

We first set up some notation that will be needed later. Let $\HHH_B:=\mathrm{Homeo}^+(S^2,B)$ be the group of orientation-preserving homeomorphisms $(S^2,B)\to(S^2,B)$ which fix $B$ pointwise, and define $\HHH_A:=\mathrm{Homeo}^+(S^2,A)$ analogously.  Let $\HHH_{\{f^{-1}(B)\}}$ denote the group of orientation-preserving homeomorphisms $(S^2, f^{-1}(B)) \to (S^2, f^{-1}(B))$ that send the set $f^{-1}(B)$ to itself, possibly via a nontrivial permutation. 

Next, we will make use of the following fact from algebra. If $A, B, C$ are groups with $A < B$, $n:=[
A:B]<\infty$, and $p: B \to C$ is a surjective homomorphism, then the image $p(A)$ has finite index in $C$.  To see this, write $B=b_1A \sqcup \ldots \sqcup b_nA$ and apply the homomorphism $p$.  We will also use the fact that the intersection of a finite collection of finite-index subgroups is again finite-index.

Let $\QQQ < \HHH_B \times \HHH_{\{f^{-1}(B)\}}$ be the subgroup of those ordered pairs of homeomorphisms $(h, \tilde{h})$ such that $h \circ f = f \circ \tilde{h}$.  Let $\mathrm{Lift}_f$ denote the image of $\QQQ$ under the projection to $\HHH_B$.

The group $\QQQ$ acts on the set $f^{-1}(B)$.  The intersection
\[ \QQQ_A:= \bigcap_{a \in A} \QQQ_a\]
of the stabilizers $\QQQ_a, a \in A$, is the subgroup of $\QQQ$ consisting of pairs $(h, \tilde{h})$ for which $h \circ f = f \circ \tilde{h}$ and for which the upstairs map $\tilde{h}$ fixes $A$ pointwise.
Thus $[\QQQ:\QQQ_A] < \infty$. Let $\mathrm{Lift}_{f,A}$ denote the image of $\QQQ_A$ under the projection to $\HHH_B$.

Note that both $\mathrm{Lift}_{f}$ and $\mathrm{Lift}_{f,A}$ are ``saturated'' with respect to homotopy; that is, for $h\in \mathrm{Lift}_{f}$, $g\in \mathrm{Lift}_{f,A}$, and for all $\iota_0\in \HHH_B$ isotopic to the identity relative to $B$, we have $\iota_0\circ h\in\mathrm{Lift}_{f}$ and $\iota_0\circ g \in \mathrm{Lift}_{f,A}$. This follows from the homotopy-lifting property for the covering map $f:(S^2,f^{-1}(B))\to (S^2,B)$. Indeed, there is a homeomorphism $\iota:(S^2,f^{-1}(B))\to (S^2,f^{-1}(B))$ which is isotopic to the identity $(S^2,f^{-1}(B))\to (S^2,f^{-1}(B))$ relative to $f^{-1}(B)$, so that $\iota_0\circ f = f\circ \iota$. Clearly, $(\iota_0\circ h,\iota\circ \tilde h)\in \QQQ$ and $(\iota_0\circ g,\iota\circ \tilde g)\in \QQQ_A$ (where $\tilde h$ and $\tilde g$ are the respective lifts of $h$ and $g$); therefore $\iota_0\circ h\in \mathrm{Lift}_{f}$ and $\iota_0\circ g\in \mathrm{Lift}_{f,A}$. 

Now consider the following commutative diagram:

\[
\xymatrix{
G_A & \ar[l]^{}  \QQQ_A \ar[d]^{*} \ar@{^{(}->}[r]       & \QQQ \ar@{^{(}->}[r]\ar[d]     & \HHH_B \times \HHH_{\{f^{-1}(B)\}} \ar[d]\\
         & \ar[d] \ar@{^{(}->}[r] \mathrm{Lift}_{f,A} & \ar[d] \ar@{^{(}->}[r] \mathrm{Lift}_f                                    & \HHH_B\ar[d] \\
         & H_f \ar[uul]^{\phi_f}        \ar@{^{(}->}[r]          & \mathrm{Lift}_f/\Homeo_0^+(S^2, B) \ar@{^{(}->}[r]     & G_B \\
}
\]
The right-pointing horizontal arrows are inclusions. The left-pointing horizontal arrow at top left  is the composition of (i) projection onto the second factor $\HHH_{\{f^{-1}(B)\}}$ (by construction, this yields an element of $\HHH_A$), with (ii) recording the isotopy class (yielding an element of $G_A$). The vertical arrows in the top row are projections onto the first  factor $\HHH_B$. The vertical arrows in the bottom row record the isotopy class as an element of $G_B$. Thus the vertical arrows are surjective by definition of their images.

Since $[\QQQ:\QQQ_A] < \infty$ we conclude using the above algebra fact that $[\mathrm{Lift}_f/\Homeo_0^+(S^2, B): H_f] < \infty$. In the remainder of this paragraph, we prove that
$[G_B:\mathrm{Lift}_f/\Homeo_0^+(S^2, B)] < \infty$, from which it follows that $[G_B: H_f]$ is finite. Indeed, let $d:= \mathrm{deg} f$, and consider all subgroups of $\pi_1(S^2,B)$ (up to conjugacy) of index $d$. Since $\pi_1(S^2,B)$ is finitely generated, there are finitely many such subgroups, hence finitely many such conjugacy classes.    The restriction $f: S^2\setminus f^{-1}(B) \to S^2\setminus B$ is a covering map, and so determines a conjugacy class $\xi$ of subgroups of index $d$ in $\pi_1(S^2, B)$ corresponding to loops in $S^2 \setminus B$ that lift to loops in $S^2 \setminus f^{-1}(B)$.   A homeomorphism $h \in \HHH_B$ determines an outer automorphism of $\pi_1(S^2, B)$, so we obtain a homomorphism $G_B \to \mathrm{Out}(\pi_1(S^2, B))$ and a corresponding action of $G_B$ on the finite set of conjugacy classes of subgroups of index $d$. The stabilizer $(G_B)_\xi$ of $\xi$ is therefore a finite index subgroup of $G_B$.  From elementary covering space theory, a sufficient (and, actually, necessary) condition for a homeomorphism $h \in \HHH_B$ to lie in the subgroup $\mathrm{Lift}_f$ is that the class $[h]$ of $h$ in $G_B$ lie in the stabilizer $(G_B)_\xi$. We conclude that $\mathrm{Lift}_f/\Homeo_0^+(S^2, B) = (G_B)_\xi$ has finite index in $G_B$.

Next, we show that the vertical arrow labelled $*$ is an isomorphism. It is surjective by definition. By definition, an element of the kernel is a pair of homeomorphisms of the form $(\id_{S^2}, \tilde{h})$ with $f = f \circ \tilde{h}$ (so that $\tilde{h}$ is a deck transformation of the covering $f: S^2\setminus f^{-1}(B) \to S^2\setminus B$) and with $\tilde{h}|_A = \id_A$. We claim that this implies $\tilde{h}=\id_{S^2}$.   Since nontrivial deck transformations have no fixed-points,  we will show that $\tilde{h}$ has a fixed-point. We do this by applying the Lefschetz fixed-point formula. To deal with compact spaces, let $U=S^2\setminus N_\epsilon(B)$ (where $N_\epsilon(B)$ is an $\epsilon$-neighborhood of $B$), and $\tilde{U}:=f^{-1}(U)$; then $\tilde{h}|_{\tilde{U}}$ is a deck transformation of the covering of compact spaces $f: \tilde{U} \to U$. Since $\#A \geq 3$ and by assumption $\tilde{h}$ fixes each element of $A$, we have that the trace of $\tilde{h}_*$ on $H_1(\tilde{U}, \Z)$ is at least $2$. By the Lefschetz fixed-point formula, this implies $\tilde{h}$ has a fixed-point. 

Next, we define the virtual endomorphism $\phi_f$. Its domain is the subgroup $H_f$. Given a class $\theta \in H_f$, to define $\phi_f(\theta) \in G_A$, choose a representing homeomorphism $h \in \mathrm{Lift}_{f,A}$. Since the vertical arrow $(*)$ is an isomorphism, there is a unique lift $\tilde{h}$ of $h$ under $f$ with the property that it fixes $A$ pointwise. We define $\phi_f(\theta):=[\tilde{h}] \in G_A$.

We now show in two ways that $\phi_f$ is well-defined.  First: given two pairs $(h_1, \tilde{h}_1), (h_2, \tilde{h}_2) \in \QQQ_A$ for which $[h_1]=[h_2]$ in $G_B$, we must show $[\tilde{h}_1]=[\tilde{h}_2] \in G_A$. Consider the product $(g, \tilde{g}):=(h_2^{-1}\circ h_1, \tilde{h}_2^{-1}\circ \tilde{h_1})$.  Since up to isotopy $f \circ \tilde{g}=f$, we  have
\[
\sigma_{f \circ \tilde{g}}=\sigma_{f}\quad \text{as maps}\quad \TTT_B\to \TTT_A.
\]
So for all $\tau\in \TTT_B$, 
\[
\sigma_{\tilde{g}}(\sigma_f(\tau))=\sigma_f(\tau)
\]
which implies that $\sigma_{\tilde{g}}$ has a fixed-point.
The action of $\sigma_{\tilde{g}}$ on $\TTT_A$ depends only on the class $[\tilde{g}] \in G_A$, and this action is via a deck transformation of the covering $\TTT_A \to \MMM_A$.  Using again the fact that a deck transformation with a fixed-point is the identity, this implies $[\tilde{g}]=\id \in G_A$.  Second: an isotopy connecting $h_1$ and $h_2$, which is constant and equal to the identity on $B$, lifts by $f$ to an isotopy which is constant on $f^{-1}(B)$, and which connects $\tilde{h}_1$ to a lift $\tilde{h}'_2$ of $h_2$. Since $\tilde{h}_1$ fixes $A\subset f^{-1}(B)$ pointwise, so does $\tilde{h}'_2$. We have shown that such a lift is unique, hence $\tilde{h}'_2=\tilde{h}_2$ and $[\tilde{h}_1]=[\tilde{h}_2] \in G_A$.

Finally, we claim that $\phi_f$ is a homomorphism. This follows from the definition of $\phi_f$ and the fact that $\mathrm{Lift}_{f,A}$ is a group under coordinatewise composition.
\qed
\gap

\noindent{\bf Restricted virtual homomorphism and other definitions.} Here, we briefly comment on the difference between the definition of virtual homomorphism given here and that of the virtual endomorphism given in \cite{kmp:tw}.  

In   \cite{kmp:tw}, the discussion treats only the dynamical setting of Thurston maps, and the corresponding virtual endomorphism is defined differently. Here is the connection. 

Given an admissible cover $f:(S^2, A) \to (S^2, B)$, the {\em restricted} virtual endomorphism, $\phi_f': G_B \dasharrow G_A$ is defined as follows.  Set $A' = f^{-1}(B)$, and let $\phi_f'': G_B \to G_{A'}$ be the virtual endomorphism given by Proposition \ref{virtual}; denote  its domain $H_f'$. The restricted virtual endomorphism $\phi_f'$ is defined as the composition 
\[ H_f' \stackrel{\phi_f''}{\longrightarrow}G_{A'} \to G_A\]
where $G_{A'} \to G_A$ is the map induced by forgetting points in $A'\setminus A$. 

Now suppose $A=B=P\supset P_f$ where $f$ is a Thurston map. The virtual endomorphism $G_P \dasharrow G_P$ of \cite{kmp:tw} coincides with the restricted virtual endomorphism $\phi_f': G_P \dasharrow G_P$. 

\gap

The virtual homomorphism changes in a predictable way under pre- and post-composition by homeomorphisms. 

\begin{lemma}\label{prepost}
Let $f:(S^2,A)\to (S^2,B)$ be an admissible covering map, let $i:(S^2,B)\to (S^2,B)$ be an orientation-preserving homeomorphism\ns{,} which maps the set $B$ to itself (not necessarily pointwise), let $j:(S^2,A)\to (S^2,A)$ be an orientation-preserving homeomorphism\ns{,} which maps the set $A$ to itself (not necessarily pointwise). Then $i\circ f\circ j:(S^2,A)\to (S^2,B)$ is an admissible cover with associated virtual homomorphism 
\[
\phi_{i\circ f\circ j}:G_B\dashrightarrow G_A\text{ with domain }H_{i\circ f\circ j},
\]
and
\[
H_{i\circ f \circ j}=i\circ H_f\circ i^{-1}\quad\text{and}\quad \phi_{i\circ f \circ j}(h) = j^{-1}\circ \phi_f(i^{-1}\circ h\circ i) \circ j, \quad h\in H_{i\circ f \circ j}.
\]
\end{lemma}
\pf This follows immediately from the definitions. \qed
\gap

\noindent{\bf Functional identities.} The following are two fundamental functional identities relating the virtual homomorphism, the linear map, and the pullback map. 

The result \cite[Thm. 1.2]{kmp:tw}, phrased and proved for Thurston maps, generalizes in a completely straightforward way to the following result; note that the more restrictive domain is referred to below. 

\begin{thm}
\label{thm:philambda} 
For any multitwist $M_w \in \dom(\phi_f')$, we have 
\begin{equation}
\label{eqn:philambda}
\phi_f'(M_w) = \phi_f(M_w) = M_{\lambda_f(w)}.
\end{equation}
\end{thm} 

From the definitions (see also  \cite{bekp} and \cite{koch:criticallyfinite}) we have 
\begin{equation}
\label{eqn:functional}
\sigma_f(h\cdot\tau) = \phi_f(h)\cdot \sigma_f(\tau), \;\; \forall \ h \in H_f.
\end{equation}
\gap

\noindent{\bf The Hurwitz space $\WWW_f$.} 
This discussion is extracted from \cite{koch:criticallyfinite}. Consider the space 
\[
\mathrm{Rat}_d\times \left(\P^1\right)^B\times\left(\P^1\right)^A
\]
where $\mathrm{Rat}_d$ denotes the space of rational maps of degree $d$, $\left(\P^1\right)^B$ denotes the space of all {\em{injective}} maps $\varphi:B\rightarrowtail\P^1$, and $\left(\P^1\right)^A$ denotes the space of all {\em{injective}} maps $\varphi:A\rightarrowtail\P^1$ (we are abusing notation). This space is a smooth affine variety. 

The group $\mathrm{Aut}(\P^1)\times\mathrm{Aut}(\P^1)$ acts on the space in the following way: 
\[
(\mu,\nu)\cdot (F,i_B,j_A)\mapsto \left(\mu\circ F\circ \nu^{-1},\mu\circ i_B,\nu\circ j_A\right).
\]
This action is free, and we consider the geometric quotient 
\[
\left(\mathrm{Rat}_d\times \left(\P^1\right)^B\times\left(\P^1\right)^A\right)/\left(\mathrm{Aut}(\P^1)\times\mathrm{Aut}(\P^1)\right);
\]
it follows from geometric invariant theory that this is a complex manifold of dimension $2d+1+\#A+\#B-6$ (for details, see \cite{koch:criticallyfinite}). 

Consider the map 
\[
\omega_f:\TTT_B\to \left(\mathrm{Rat}_d\times \left(\P^1\right)^B\times\left(\P^1\right)^A\right)/\left(\mathrm{Aut}(\P^1)\times\mathrm{Aut}(\P^1)\right)
\]
given by $\omega_f:[\phi]\mapsto \left[F,\phi|_B,\psi|_A\right]$ where 
\begin{itemize}
\item $\phi:(S^2,B)\to (\P^1,\phi(B))$ and $\psi:(S^2,A)\to (\P^1,\psi(A))$  are orientation-preserving homeomorphisms, and 
\item $F:=\phi \circ f \circ \psi^{-1}:(\P^1,\psi(A))\to (\P^1,\phi(B))$ is the rational map at right in the diagram defining the pullback map $\sigma_f$.
\end{itemize}
In other words, $\omega_f$ records the algebraic data of $\sigma_f$, up to appropriate equivalence. 

Following \cite{koch:criticallyfinite}, the {\em{Hurwitz space}} associated to $f:(S^2,A)\to (S^2,B)$ is $\WWW_f:=\omega_f(\TTT_B)$. It is a complex manifold of dimension $\#B-3$, and the map $\omega_f:\TTT_B\to \WWW_f$ is a holomorphic covering map. 
\gap

\begin{figure}[h]

\[
\xymatrix{  & \TTT_B \ar[dd]_{\pi_B}\ar[rr]^{\sigma_f} \ar[dr]^{\omega_f} & &
\TTT_A \ar[dd]^{\pi_A} \\
&&\WWW_f\ar[dl]_Y \ar[dr]^X &\\
& \MMM_B & & \MMM_A}
\]
\caption{\sf The fundamental diagram associated to any admissible covering map $f:(S^2,A)\to (S^2,B)$; this diagram commutes.}\label{fund}
\label{figure:w} 
\end{figure}
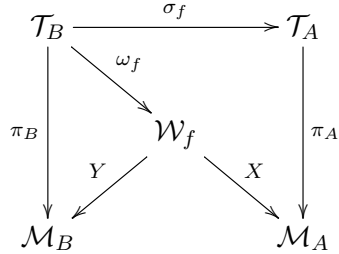

\noindent{\bf The moduli space correspondence.} The Hurwitz space $\WWW_f$ associated to $f:(S^2,A)\to (S^2,B)$ is equipped with two holomorphic maps (see Figure \ref{figure:w}): 
\[
X:\WWW_f\to \MMM_A\quad\text{given by}\quad [F,\phi|_B,\psi|_A]\mapsto [\psi|_A]\quad \text{and}
\]
\[
Y:\WWW_f\to \MMM_B\quad\text{given by}\quad  [F,\phi|_B,\psi|_A]\mapsto [\phi|_B].
\]

We may think of these as comprising one multivalued map $X \circ Y^{-1}: \MMM_B \rightrightarrows \MMM_A$; note that the direction of this correspondence is the same as that of $\sigma_f:\TTT_B\to\TTT_A$, i.e. it corresponds to pulling-back complex structures.  We call the pair of maps 
\[
\left(X:\WWW_f\to \MMM_A,\quad Y:\WWW_f\to \MMM_B\right)
\]
the {\em moduli space correspondence} associated with $f:(S^2,A)\to (S^2,B)$. The following results are proved in \cite{koch:criticallyfinite}.
\begin{prop}\label{Wquotient}
The map $\pi_B$ factors as $\pi_B=Y\circ \omega_f$, and the space $\WWW_f$ is isomorphic (as a complex manifold) to $\TTT_B/H_f$. 
\end{prop}
\begin{cor}\label{Wfinitecover}
The map $Y:\WWW_f\to \MMM_B$ is a finite covering map. 
\end{cor}
\pf This follows from Proposition \ref{virtual}. 
\qed 
\gap

On the one hand, the map $Y: \WWW \to \MMM_B$ is always a covering map.  On the other hand, the map $X:\WWW_f\to \MMM_A$ may have diverse properties: it can be e.g. constant, injective but not surjective, or a surjective ramified covering; see \cite{bekp}.
\gap
 
We say that admissible covers $f:(S^2,A)\to (S^2,B)$ and $g:(S^2,A)\to (S^2,B)$ are $(A,B)$-{\em Hurwitz equivalent} if there are homeomorphisms $h\in\HHH_B$, and $\tilde h\in\HHH_A$ so that $h\circ f=g\circ \tilde h$.  
\gap

While $\sigma_f$ depends on the admissible cover $f:(S^2,A)\to (S^2,B)$, the correspondence on moduli space depends only on the $(A,B)$-Hurwitz class of $f$. 
\begin{prop}
The Hurwitz space $\WWW_f$ is equal to the Hurwitz space $\WWW_g$ if and only if $f$ and $g$ are $(A,B)$-Hurwitz equivalent. 
\end{prop}
\gap

\noindent{\bf Remark:} The classical notion of Hurwitz equivalence requires only that $h\circ f=g\circ \tilde h$ for some pair of orientation-preserving homeomorphisms $h, \tilde{h}$ which are not required to fix any finite set pointwise. Thus, $(A,B)$-Hurwitz equivalence is a finer equivalence relation than classical Hurwitz equivalence.  For example, any pair of quadratic maps $f, g$ are Hurwitz equivalent, but it is possible to arrange a suitable choice of $A=B=P_f = P_g$ so that $f, g$ are not $(A,B)$-Hurwitz equivalent; see \cite[Remark 7.3.3]{koch:thesis} and \cite[Example 2.9]{koch:criticallyfinite}.
\gap

\noindent{\bf Two definitions of the virtual homomorphism.}  Choose a basepoint $w_\circledast \in \WWW$ and set $\mu_\circledast:=X(w_\circledast), \tilde{\mu}_\circledast :=Y(w_\circledast)$. Using the identity maps to define corresponding basepoints for $\TTT_A, \TTT_B$, we obtain identifications (see \cite[Section 2.2]{lodge:thesis})
$G_A \leftrightarrow \pi_1(\MMM_A, \tilde{\mu}_\circledast)$ and $G_B \leftrightarrow \pi_1(\MMM_B, \mu_\circledast)$. 
It can be shown (cf. \cite[Theorem 2.6]{lodge:thesis}) that the virtual homomorphism $\phi_f: G_B \dasharrow G_A$, under these identifications,  coincides with the induced maps $X_* \circ (Y_*)^{-1}$ on fundamental groups. 
 
 \gap
 
 \noindent{\bf Augmented Teichm\"uller space and the extension of $\sigma_f$.}  The proof of \cite[Prop. 4.3]{selinger:boundary} shows that $\sigma_f$ is uniformly Lipschitz with respect  to the WP metric, and therefore extends to the WP completion. The fundamental diagram becomes 
\begin{equation}
\label{eqn:cptfunddiag} 
\xymatrix{  & \cl{\TTT}_B \ar[dd]_{\pi_B}\ar[rr]^{\sigma_f} \ar[dr]^{\omega_f} & &
\cl{\TTT}_A \ar[dd]^{\pi_A} \\
&&\cl{\WWW}_f\ar[dl]_Y \ar[dr]^X &\\
& \cl{\MMM}_B & & \cl{\MMM}_A}
\end{equation}
where now $ \cl{\MMM}_B$, $\cl{\MMM}_A$ are the {\em Deligne-Mumford compactifications}.

The proof of  \cite[Prop. 6.1]{selinger:boundary} shows 
\begin{equation}
\label{eqn:stratamap}
\sigma_f(\TTT_B^\Gamma)\subset \TTT_A^{f^{-1}(\Gamma)}.
\end{equation}

\subsection{Dynamics}\label{dynamics} We briefly mention the relevant objects in the dynamical setting. 
\gap

\noindent{\bf Thurston maps.} We will often be concerned with the special case where $A=B=P \supset P_f$.
We will require the notion of the orbifold associated to a Thurston map $f:(S^2,P)\to (S^2,P)$. Following \cite{DH1}, we note that there is a smallest function $\nu_f$ over all functions $\nu:S^2\to\N\cup\{\infty\}$ such that 
\begin{itemize}
\item $\nu(x)=1$ if $x\notin P$, and 
\item $\nu(x)$ is a multiple of $\nu(y) \cdot \deg(f,y)$ for all $y\in f^{-1}(x)$; here the second factor is the local degree of $f$ at $y$.
\end{itemize}
Indeed, one simply sets $\nu(x)$ to be the least common multiple of the set of local degrees of iterates of $f$ at all points that are iterated preimages of $x$. 
The {\em{orbifold}} of $f$ is $O_f(S^2,\nu_f)$; it is {\em{hyperbolic}} if the Euler characteristic
\[
\chi(O_f)=2-\sum_{x\in P}\left(1-\frac{1}{\nu_f(x)}\right)
\]
is negative. We say that $f$ has {\em{Euclidean orbifold}} if it does not have hyperbolic orbifold. If $\#P_f\geq 5$, then the orbifold of $f$ is necessarily hyperbolic. Maps with Euclidean orbifolds are classified in \cite{DH1}. 
\gap

\noindent{\bf Pullback relation on curves.} Just as in the nondynamical setting, the pullback relation induces a {\em pullback function} $f^{-1}: \MMM\SSS_P \to \MMM\SSS_P$.  A multicurve is {\em invariant} if $f^{-1}(\Gamma) \subset \Gamma$ or $f^{-1}(\Gamma) = \emptyset$; it is {\em completely invariant} if $f^{-1}(\Gamma) = \Gamma$.  

\gap 

\noindent{\bf The Thurston linear transformation.} The {\em Thurston linear transformation} 
\[  \lambda_f: \R[\SSS_P] \to \R[\SSS_P]\]
is defined just as before, with $A=B=P$. 

There is one key difference in the dynamical setting, which is 
{\em Thurston's Characterization Theorem} \cite[Theorem 1]{DH1}.  

\begin{thm}[Thurston's characterization]
\label{thm:Thurston_characterization}
If $\mathcal{O}_f$ is hyperbolic, then $f$ is equivalent to a rational map $R$ if and only if for every  invariant (equivalently, every) multicurve $\Gamma$, the spectrum of the linear map 
\[
\lambda_{f,\Gamma}:=\lambda_f|\R[\Gamma]: \R[\Gamma]\to\R[\Gamma]
\]
 lies in the open unit disk; in this case, $R$ is unique, up to conjugation by M\"obius transformations.   
\end{thm}
If there is an invariant multicurve $\Gamma\subset S^2\setminus P$ for which the criterion in Theorem \ref{thm:Thurston_characterization} fails, then $\Gamma$ is a called a {\em{Thurston obstruction}}. 

\gap

\noindent{\bf The pullback map.} The pullback map $\sigma_f:\TTT_P\to\TTT_P$ is defined as before.  

When $\OOO_f$ is hyperbolic, the criterion in Theorem \ref{thm:Thurston_characterization} is equivalent to the existence of a fixed point of $\sigma_f$ in $\TTT_P$.  Some $k$th iterate $\sigma^k_f$ is a strict, non-uniform contraction (when $P=P_f$, \cite[Prop. 3.3(b)] {DH1} asserts that $k=2$ will do).  So if a fixed point exists, then it is unique, $f$ is equivalent to a unique rational map $R$, and the projection of the fixed-point to moduli space corresponds to the geometry of the corresponding subset of the dynamical plane of $R$. When $\OOO_f$ is Euclidean, the relationship between the existence of Thurston obstructions and the dynamics on Teichm\"uller space is more subtle; cf. \cite{selinger:euclidean}. However, with the exception of the well-known {\em integral Latt\'es maps} induced by the endomorphisms $z \mapsto n\cdot z$ on complex tori, a fixed-point of $\sigma_f$, if it exists, is unique.  
\gap

\noindent{\bf The virtual endomorphism.}  There are a finite index subgroup $H_f<G_P$ and a virtual endomorphism 
\[
\phi_f:G_P\dashrightarrow G_P\quad\text{given by}\quad \phi_f:[h]\mapsto [h']
\]
where $h,h'\in \Homeo^+(S^2, P)$, and $h \circ f = f\circ h'$. 
\gap

\noindent{\bf The Hurwitz space $\WWW_f$.} We have the same fundamental diagram as in Figure \ref{fund}, with $A=B=P$, and we have maps 
\[
X:\WWW_f\to \MMM_P\quad\text{and}\quad Y:\WWW_f\to \MMM_P.
\]
In the dynamical setting, we can iterate the moduli space correspondence $X\circ Y^{-1}:\MMM_P\rightrightarrows\MMM_P$.
\gap

\noindent{\bf Augmented Teichm\"uller space and the extension of $\sigma_f$.} Considering the Diagram (\ref{eqn:cptfunddiag}), the proof of  \cite[Prop. 6.1]{selinger:boundary} shows again that 
$\sigma_f(\TTT_P^\Gamma)\subset \TTT_P^{f^{-1}(\Gamma)}$.  In particular, completely invariant multicurves (those satisfying $f^{-1}(\Gamma) = \Gamma$) correspond to strata invariant under $\sigma_f$.    Selinger \cite[Thm. 10.4]{selinger:boundary} shows that if $f$ is obstructed, then under suitable hypotheses (all pieces in the corresponding canonical decomposition having hyperbolic orbifold), again there is a unique fixed-point $\hat{\tau} \in \cl{\TTT}_P$ to which all iterates converge. However, there exist both obstructed and unobstructed maps $f$ for which $\sigma_f$ has periodic points in $\bdry \TTT_P$. Indeed, understanding the dynamical behavior of $\sigma_f: \cl{\TTT}_P \to \cl{\TTT}_P$ is a main motivation for this work. 
\gap

\section{When $\sigma_f:\bdry\TTT_B \to \bdry\TTT_A$}

The main result of this section is the following characterization of when the pullback map preserves the WP boundary.

\begin{thm}
\label{thm:bdry}
Let $f:(S^2,A)\to (S^2,B)$ be an admissible cover. The following are equivalent.
\be
\item  For each $\Gamma \in  \MMM\SSS_B$, $f^{-1}(\Gamma)\neq \emptyset$. 
\item For each $\gamma \in \SSS_B$, $\gamma\notin \ker(\lambda_f)$
\item $\ker(\phi_f) \intersect \Tw^+\subset G_B$ is trivial. 
\item  The pullback map $\sigma_f$ on Teichm\"uller space satisfies $\sigma_f: \bdry\TTT_B \to \bdry\TTT_A$.
\item  The pullback correspondence $X\circ Y^{-1}:\MMM_B\rightrightarrows\MMM_A$ on moduli space is proper.
\eb
If $\#A=\#B$, and if any of the conditions above hold, the pullback correspondence $X\circ Y^{-1}:\MMM_B\rightrightarrows\MMM_A$ is surjective.
\end{thm}

In (5), properness means that for each compact $K \subset \MMM_A$, the set $Y(X^{-1}(K)) \subset \MMM_B$ is compact. Since $Y$ is a finite covering map, this is equivalent to the properness of the map $X$.  
\gap

\pf The equivalence of $(1)$ and $(2)$ is immediate from the definitions. 

The equivalence of $(2)$ and $(3)$ follows immediately from the definitions, Equation \ref{eqn:philambda}, the nonnegativity of $\lambda_f$, and the fact that $\phi_f$ preserves the positivity of twists. Indeed, from \cite[Thm. 1.2]{kmp:tw} we have the following. Suppose $\Gamma=\{\gamma_1, \ldots, \gamma_n\}$ is a nonempty multicurve on $S^2\setminus B$, $a_1, \ldots, a_n \in \Z_{>0}$, $w:=\sum a_i\gamma_i \in \Z_+[\SSS_B]$, and $M_w:=T_1^{a_1}\cdot \ldots \cdot T_n^{a_n}$ is the corresponding positive multitwist; here $T_i$ is the right Dehn twist about $\gamma_i$. If $M_w \in H_f'$, then by Equation (\ref{eqn:philambda}) 
\[ \phi_f(M_w) = M_{\lambda_f(w)}.\]
So 
\[ \phi_f(M_w) = \id \iff \lambda_f(w)=0.\]
Since $\lambda_f$ is a nonnegative linear operator and $a_i > 0$ for each $i$, we conclude from the definition of $\lambda_f$ that 
\[ \phi_f(M_w) = \id \iff f^{-1}(\Gamma) = o.\]

%

We now show (1) $\iff$ (4). Recall that the WP boundary $\bdry \TTT_B$ is a union of strata corresponding to nonempty multicurves on $S^2\setminus B$. By Equation (\ref{eqn:stratamap}) 
\[ \sigma_f(\TTT_B^\Gamma) \subset \TTT_A^{f^{-1}(\Gamma)}.\]
In particular, 
\[ \sigma_f(\TTT_B^\Gamma) \subset \bdry \TTT_A \iff f^{-1}(\Gamma) \neq o,\]
which yields $(1) \iff (4)$.

Now we show that failure of $(4)$ implies failure of $(5)$.  Suppose for some nonempty multicurve $\Gamma$ there exist $\tau_n \to \overline{\tau} \in \TTT_B^\Gamma$ while $\tilde{\tau}_n:=\sigma_f(\tau_n) \to \tilde{\tau}:=\sigma_f(\overline{\tau}) \in \TTT_A$. But then in moduli space $\mu_n:=\pi_B(\tau_n) \to \infty $ while $\tilde{\mu}_n:=\pi_A(\sigma_f(\tau_n)) \subset X\circ Y^{-1}(\mu_n)$ does not, showing that properness of the pullback correspondence fails. 


Now suppose that $(5)$ fails. Then there exist $\tilde{\mu}_n \to \tilde{\mu}$ in $\MMM_A$ and and $\mu_n \to \overline{\mu} \in \partial \MMM_B$ such that $\tilde{\mu}_n \in X\circ Y^{-1}(\mu_n)$.  Consider Diagram~\ref{eqn:cptfunddiag}.  There exists $\cl{w} \in \overline{\WWW_f}$ with $Y(\cl{w})=\overline{\mu}$, $X(w)=\tilde{\mu}$.  We see from Diagram~\ref{eqn:cptfunddiag} that for any $\cl{\tau} \in \partial\TTT_B$ with $\omega_f(\cl{\tau})=\cl{w}$, $\sigma_f(\cl{\tau}) \in \pi_A^{-1}(\tilde{\mu}) \in \TTT_A$ and $(4)$ fails.

It remains to prove surjectivity of $X\circ Y^{-1}$ under the assumptions $\# A=\# B$, and that e.g. Condition $(5)$ holds. Since $X\circ Y^{-1}$ is proper, it is not constant. 

Recall that $\MMM_B$ is isomorphic to a hyperplane complement; in particular, it is a Stein manifold and, hence, the finite cover $\WWW_f$ of $\MMM_B$ also is a Stein manifold of the same dimension. Therefore, the preimage under $X$ of any point is totally disconnected and, by properness, compact. Since $X$ is analytic, we conclude that the preimage is discrete as well, hence finite. Because $\# A=\# B$, the complex manifolds $\WWW_f$ and $\MMM_A$ have equal dimension, so the map $X:\WWW_f\to \MMM_A$ is open because it has discrete fibers. A continuous proper map between locally compact Hausdorff spaces is closed.  Hence the image $X(\WWW_f)$ is both open and closed in $\MMM_A$ and so $X:\WWW_f\to\MMM_A$ is surjective.
\qed

\section{When $\sigma_f:\TTT_B\to\TTT_A$ is constant} 

\begin{thm}\label{thm:constant}
Let $f:(S^2,A)\to (S^2,B)$ be an admissible cover. The following are equivalent:
\begin{enumerate}
\item $\pullback$ is constant
\item $\lambda_f:\R[\SSS_B]\to\R[\SSS_A]$ is constant
\item $\phi_f:G_B\dashrightarrow G_A$ is constant
\item $\sigma_f:\TTT_B\to\TTT_A$ is constant
\item $X\circ Y^{-1}:\MMM_B\rightrightarrows\MMM_A$ is constant
\end{enumerate}
\end{thm}

In \cite{bekp}, the previous theorem is proved in the dynamical setting for $f$ a Thurston map. The same ideas for the proof apply in the nondynamical case. However, in \cite{bekp}, there is a mistake in the proof that $(2)\implies (3)$.  The assumption $(2)$ is equivalent to the assumption that every curve, when lifted under $f$, becomes inessential or peripheral.  Even if this holds, it need not be the case that every Dehn twist lifts under $f$ to a pure mapping class element.  We give an explicit example after the proof of Theorem \ref{thm:constant}. 

\pf 

In \cite{bekp} the logic was: $(1)\implies(2)\implies(3)\implies(4)$, and  failure of $(1)$ implies failure of $(4)$. 
Also, condition $(1)$ is stated as ``$\pullback$ is empty''; and condition  $(5)$ is omitted. 

Here is the revised logic: 
\begin{itemize}
\item $(4)\iff (5)$,
\item $(1)\iff (2)$, 
\item $(3)\implies (2)$,
\item $(3)\iff (4)$, and 
\item failure of $(4)$ implies failure of $(1)$. 
\end{itemize}

We immediately have that $(4)\iff (5)$ from the fundamental diagram, Figure \ref{fund}, and $(1)\iff (2)$ follows immediately from the definitions.  

To show that ${(3)\implies (2)}$, we show failure of $(2)$ implies failure of $(3)$. If $\lambda_f$ is not constant, then there exists a simple closed curve $\gamma\in\SSS_B$ which has an essential, nonperipheral  simple closed curve $\delta\in\SSS_A$ as a preimage under $f$.  Some power of the Dehn twist about $\gamma$ lifts under $f$ to a product of nontrivial Dehn twists.  The hypothesis implies that the lifted map is homotopically nontrivial, so $\phi_f$ is not constant. 

For the remaining implications, we will make use of the following facts: First recall the functional identity from (\ref{eqn:functional}):
\[
\sigma_f(h\cdot\tau) = \phi_f(h)\cdot \sigma_f(\tau), \;\; \forall \ h \in H_f,
\]
recall that $W_f$ is isomorphic (as a complex manifold) to $\TTT_B/H_f$ (Proposition \ref{Wquotient}), and recall Corollary \ref{Wfinitecover} which states that $Y:\WWW_f\to \MMM_B$ is a finite cover.

Second, note that a bounded holomorphic function on a finite cover of $\MMM_B$ is constant.  To see this, recall that $\MMM_B$ is isomorphic to the complement of a finite set of hyperplanes in $\C^n$ where $n=\#B-3$.  Let $L\in\C^n$ be any complex line which intersects $\MMM_B$.  This intersection  is isomorphic to a compact Riemann surface punctured at finitely many points.  If $\widetilde{L}$ is any component of the preimage of $L$ under a finite covering, then $\widetilde{L}$ is also isomorphic to a compact Riemann surface punctured at finitely many points.   By Liouiville's theorem, the function is constant on $\widetilde{L}$. Since $L$ is arbitrary, the function is locally constant, hence constant.  

To establish ${(3)\implies(4)}$, suppose $(3)$ holds.  Then  $\sigma_f: \TTT_B\to \TTT_A$ descends
to a holomorphic map 
\[
\overline{\sigma}_f: \WWW_f \to \TTT_A.
\]

But it is well-known (see e.g. \cite{hubbard:teich_vol1}) that $\TTT_A$ is isomorphic as a complex manifold to a bounded domain of $\C^n$, so the discussion above implies that $\sigma_f$ is constant.  

To establish that ${(4)\implies (3)}$, suppose $h \in H_f$.   If $\sigma_f \equiv \tau$ is constant, the deck transformation defined by $\phi_f(h)$ fixes the point $\tau$, hence must be the identity, so $\phi_f$ is constant.  

To establish that $\mathrm{not}(4) \implies \mathrm{not}(1)$, we first prove a lemma of perhaps independent interest. Below, by a bounded set we mean a set with compact closure. 

\begin{lemma}\label{unbounded} Let $f:(S^2,A)\to (S^2,B)$ be an admissible covering map. Then the image of $\sigma_f:\TTT_B\to\TTT_A$ is either a point, or unbounded in $\TTT_A$.  In the latter case, actually $\pi_A(\sigma_f(\TTT_B))$ is unbounded in $\MMM_A$. 
\end{lemma} 
\pf
The fundamental diagram in Figure \ref{fund} implies that
\[
\sigma_f(\TTT_B)\text{ is bounded in }\TTT_A\;\iff\;\pi_A(\sigma_f(\TTT_B))\text{ is bounded in }\MMM_A
\]
\[
\iff\; X(\WWW_f)\text{ is bounded in }\MMM_A.
\]
Because a bounded holomorphic function on $\WWW_f$ is constant, $X(\WWW_f)$ is bounded in $\MMM_A$ if and only if it is a single point, or (by  Figure \ref{fund}) if and only if the image of $\sigma_f$ is a single point. 
\qed

\bigskip 

Suppose now that $\sigma_f:\TTT_B\to\TTT_A$ is not constant (i.e., failure of $(4)$). Lemma \ref{unbounded} implies that $X(\WWW_f)$ is not contained in any compact subset of $\MMM_A$.   Consider the diagram in Equation (\ref{eqn:cptfunddiag}).  This means that the image $X(\overline{\WWW_f})$ contains at least one point on the boundary $\partial \MMM_A$, hence $\sigma_f(\overline{\TTT_B})$ contains at least one point in $\partial \TTT_A$. We see that there exists a boundary stratum $\TTT_B^\Gamma$ that is mapped to a boundary stratum of $\TTT_A$. Hence, $f^{-1}(\Gamma) \neq o$, which means $\pullback$ is not constant.
%
\qed

\noindent{\bf Remark:} One may prove the last implication in a more direct and elementary way by observing that if $X(\WWW_f)$ is unbounded, then there are points $\tilde{\tau}=\sigma_f(\tau)$ on whose underlying Riemann surfaces $X_{\tilde{\tau}}$ exist annuli of arbitrarily large modulus. Pushing such an annulus forward under the natural map $X_{\tilde{\tau}} \to X_\tau$ yields a fat annulus on $X_\tau$. The core curves of these annuli must be nontrivial and provide a pair $\gamma \pullback \tilde{\gamma}$. 
\medskip

\medskip

\noindent {\bf Example.}
Consider the postcritically finite rational map $f: \P^1\to\P^1$ defined by 
\[
f(z)=2i\left(z^2-\frac{1+i}{2}\right)^2,
\]
with postcritical set $P=\{0,1,-1,\infty\}$.   It factors as $f=g\circ s$ where $s(z)=z^2$. Its mapping properties are shown in Figure \ref{fig:ex1}
In this case, $\sigma_f:\TTT_P\to\TTT_P$ is constant (as proved in \cite{bekp}), $Y:\WWW_f\to\MMM_P$ is a degree $2$ covering map, and $X:\WWW_f\to\MMM_P$ is constant. Because $Y$ has degree $2$, the subgroup $H_f$ has index $2$ in $G_P$, so there is a Dehn twist $h\in \mathrm{Homeo}^+(S^2,B)$ that does not lift to a homeomorphism $h'\in \mathrm{Homeo}^+(S^2,A)$. Thus even if $\sigma_f:\TTT_B\to\TTT_A$ is constant, the index $[G_B:H_f]$ may be greater than $1$. 

One may see this directly. Let $\gamma_0$ be the boundary of a small regular neighborhood $D$ of the segment $[0,1] \subset \C$.  Let $h_0: \P^1 \to \P^1$ be the right Dehn twist about $\gamma_0$.  
\gap

\noindent{\bf Claim:} If $h_1: \P^1 \to \P^1$ satisfies $h_0 \circ f = f \circ h_1$ (i.e. $h_1$ is a lift of $h_0$ under $f$) then $h_1 \not\in \pmcg$. 
\gap

\pf We argue by contradiction.  
\begin{figure}[htbp]
\centerline{\scalebox{0.75}{\includegraphics{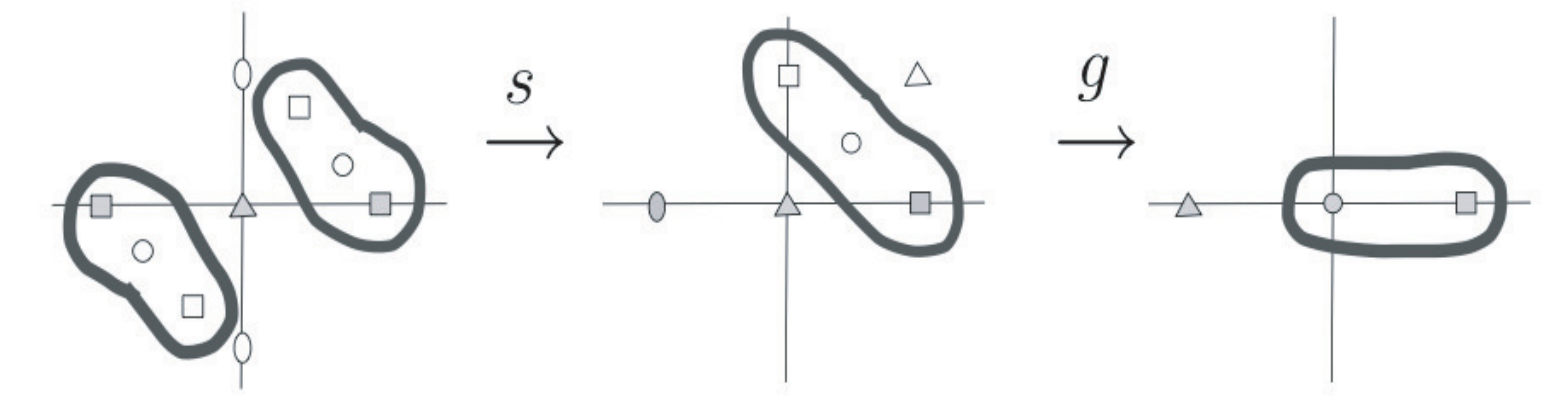}}}
\caption{The mapping properties of $f$.  The points in grey are $-1, 0, +1$. The annulus $A_0$ is shown in dark gray at right.}
\end{figure}
We may assume $h_0$ is supported on an annulus $A_0$ (see Figure 3) surrounding a bounded Jordan domain $D_0$ whose boundary is $\gamma_0$, and an unbounded region $U_0$.  Easy calculations show that the inverse image of $D_0$ under $f$ consists of two bounded Jordan domains $D_1^\pm$ each mapping as a quadratic branched cover onto $D_0$ and ramified at the points $c_\pm:=\pm \sqrt{\frac{1+i}{2}}$ (the positive sign corresponding to the root with positive real part), both of which map to the origin under $f$.  The domain $D_1^+$ contains two preimages of the point $1$, namely $+1$ and $+\frac{1+i}{\sqrt{2}}$, while its twin $D_1^-$ also contains two preimages of the point $1$, namely $-1$ and $-\frac{1+i}{\sqrt{2}}$.  The points $\pm 1 \in D_1^\pm$ belong to $P_f$, so if $h_1 \in \pmcg$ is a lift of $h_0$, then $h_1(1)=1$ and $h_1(-1)=-1$.   

Since $f: D_1^\pm - \{c_\pm\} \to D_0 - \{0\}$ are both unramified coverings, and $h_0: D_0-\{0\} \to D_0-\{0\}$ is the identity map, we conclude $h_1: D_1^\pm - \{c_\pm\} \to D_1^\pm - \{c_\pm\}$ is a deck transformation of this covering fixing a point, hence is the identity on $D_1^\pm$.

The preimage of the annulus $A_0$ is a pair of disjoint, non-nested  annuli $A_1^\pm$ with an inner boundary component $\gamma_1^\pm$ equal to $\bdry D_1^\pm$.  Since $f: A_1^\pm \to A_0$ is quadratic and unramified, and, by the previous paragraph, the restriction $h_1|_{D_1^\pm}= \id_{\gamma_1^\pm}$, we must have $h_1 \neq \id$ on the outer boundary components of $A_1^\pm$; indeed, $h_1$ there effects a half-twist.

The preimage of $U_0$ under $f$ is a single unbounded region $U_1$, which is homeomorphic to the plane minus two disks and three points; it maps in a four-to-one fashion, ramified only at the origin.  The restriction $f: U_1 - \{f^{-1}(0)\} \to U_0-\{-1\}$ is an unramified covering map, so $h_1: U_1 - \{f^{-1}(-1)\} \to U_1 - \{f^{-1}(-1)\}$ is a deck transformation of this covering.  By the previous paragraph, it is distinct from the identity.

We will obtain a contradiction by proving that $h_1: U_1 - \{f^{-1}(-1)\} \to U_1 - \{f^{-1}(-1)\}$ has a fixed point; this is impossible for deck transformations other than the identity.  We use the Lefschetz fixed point formula.  By removing a neighborhood of $\infty$ and of $-1$, and lifting these neighborhoods, we place ourselves in the setting of compact planar surfaces with boundary, so that this theorem will apply.  Under $h_1$, the boundary component near infinity is sent to itself, as are the outer boundaries of $A_1^\pm$ and the boundary component surrounding the origin (since the origin is the uniquely ramified point of $f$ over $U_0$).  The remaining pair of boundary components are permuted amongst themselves.  The action of $h_1: U_1 - \{f^{-1}(-1)\} \to U_1 - \{f^{-1}(-1)\}$ on rational homology has trace equal to either $3$ or $5$.  A fixed point thus exists, and the proof is complete.  
\qed

{\bf Remark:}  There exists a lift $h_1$ of $h_0$ under $f$.  First, there is a lift $h'$ of $h_0$ under $g$, obtained by setting $h' = \id$ on the preimage of $U_0$.  This extends to a half-twist on the preimage $A_0'$  of $A_0$ under $g$, which then in turn extends to a homeomorphism fixing the preimage $D_0'$ of $D_0$ under $g$; inside $D_0'$, this homeomorphism interchanges the points $1, i$ which are the primages of $1$.  It is then straightforward to show that $h'$ lifts under $s$ by setting $h_1=\id$ on $U_1$ and extending similarly over the annuli $A_1^\pm$ and the domains $D_1^\pm$.

\section{Noninjectivity of the virtual homomorphism} 

In this section, we begin with a nondynamical discussion about the injectivity of the virtual homomorphism $\phi_f:G_B\dashrightarrow G_A$ associated to an admissible covering map $f:(S^2,A)\to (S^2,B)$. 

It follows from Lemma \ref{prepost} that $\phi_f:G_B\dashrightarrow G_A$ is injective if and only if $\phi_{i\circ f\circ j}:G_B\dashrightarrow G_A$ is injective. 

\begin{lemma}\label{phi:sigma:inj}
Let $f:(S^2,A)\to (S^2,B)$ be an admissible cover. If $\sigma_f:\TTT_B\to\TTT_A$ is injective, then $\phi_f:G_B\dashrightarrow G_A$ is injective. 
\end{lemma}
\pf Suppose $h \in \ker(\phi_f)$. Then from the functional identity (\ref{eqn:functional}) we have $\sigma_f(h.\tau) = \sigma_f(\tau)$ for all $\tau \in \TTT_B$ and so non-injectivity of $\phi_f$ implies non-injectivity of $\sigma_f$. 
\qed
\gap

\noindent{\bf Remark.} We are not aware of any examples of an admissible covering map $f:(S^2,A)\to (S^2,B)$ where $\phi_f:G_B\dashrightarrow G_A$ is injective but $\sigma_f:\TTT_B\to \TTT_A$ is not. 

Our discussion naturally breaks up into three cases depending on $\#A$ and $\#B$. We treat the first two cases in the proposition below. The remaining case where $\#A=\#B$ is treated in Theorem \ref{thm:noninjective}.
\begin{prop}\label{BvsA}
Let $f:(S^2,A)\to (S^2,B)$ be an admissible covering map with virtual homomorphism $\phi_f:G_B\dashrightarrow G_A$. 
\begin{enumerate}
\item If $\#B>\#A$, then $\phi_f$ is not injective. 
\item If $\#B<\#A$, then $\phi_f$ is injective if $A=f^{-1}(B)$. 
\end{enumerate}
\end{prop}

\pf We begin with a proof of $(1)$ due to D. Margalit. Let $G_{\text{max}}\leq G_B$ be a maximal abelian subgroup of $G_B$. The rank of $G_{\text{max}}$ is $\#B-3$ by Theorem A in \cite{blm:duke:1983}. Let $H_{\text{max}}$ be a maximal abelian subgroup of $H_f$, the domain of $\phi_f$. Because $H_f$ has finite index in $G_B$, the rank of $H_{\text{max}}$ is also equal to $\#B-3$ (Theorem 6.4C, \cite{ivanov:mcg:1998}). If $\phi_f$ were injective, then $\phi_f(H_{\text{max}})$ would be an abelian subgroup in $G_A$ of rank $\#B-3$. But the rank of any abelian subgroup of $G_A$ is bounded above by $\#A-3$, and since $\#B>\#A$ by hypothesis, the virtual homomorphism $\phi_f$ cannot be injective. 

We now prove $(2)$. The Teichm\"uller spaces $\TTT_A, \TTT_B$ are each canonically isomorphic  to the Teichm\"uller spaces where the points in $A$ and $B$ represent punctures. Under the hypothesis that $A=f^{-1}(B)$, the pullback map $\sigma_f: \TTT_B \to \TTT_A$ is induced by lifting  complex structures under an unramified covering map. It is well-known that in this case $\sigma_f$ is a global isometry with respect to the Teichm\"uller metrics (a Teichm\"uller mapping $\psi$ corresponding to a quadratic differential $q$ lifts to a Teichm\"uller mapping corresponding to a lifted quadratic differential) and is therefore injective. By Lemma \ref{phi:sigma:inj}, $\phi_f$ is injective. 
\qed
\gap

\noindent{\bf Remark:} Condition $(2)$ is sufficient but not necessary. For example, in the dynamical setting, if $A=B=P_f$, and $f:(S^2,P_f)\to (S^2,P_f)$ is a Thurston map, and if $f$ has Euclidean orbifold, then $\phi_f$ is injective. This is immediately clear in the case where $\#P_f=3$. In the case where $\#P_f=4$, the result follows from Lemma \ref{phi:sigma:inj}, and from the fact that $\sigma_f:\TTT_{P_f}\to \TTT_{P_f}$ is an automorphism \cite{DH1}. 
\gap 

For the remaining case, we will require a nondynamical version of the notion of hyperbolic orbifold (see Section \ref{dynamics}). Let $f:(S^2,A)\to (S^2,B)$ be an admissible cover with $\#A=\#B$. We say that $f$ has {\em{potentially hyperbolic orbifold}} if there are orientation-preserving homeomorphisms $i, j$ with 
\[
j^{-1}:(S^2,A)\to (S^2,j^{-1}(A)),\quad\text{and}\quad i:(S^2,B)\to (S^2,i(B)), 
\]
so that $P:=j^{-1}(A)=i(B)$, and $i\circ f\circ j:(S^2,P)\to (S^2,P)$ is a Thurston map with hyperbolic orbifold. Note that if $\#A=\#B$, then it is easy to find homeomorphisms $j$ and $i$ so that $f:(S^2,P)\to (S^2,P)$ is a Thurston map; the content of this definition is that some such Thurston map has hyperbolic orbifold. 
\gap

\begin{thm}
\label{thm:noninjective}
Let $f:(S^2,A)\to (S^2,B)$ be an admissible cover with $\#A=\#B \geq 4$, and suppose that $f$ has potentially hyperbolic orbifold. Then the virtual homomorphism $\phi_f:G_B\dashrightarrow G_A$ is not injective: its kernel contains a pseudo-Anosov element. 
\end{thm}

Because $f$ has potentially hyperbolic orbifold, there are homeomorphisms $i, j$ with $j^{-1}:(S^2,B)\to (S^2,P)$ and $i:(S^2,A)\to (S^2,P)$ so that $F:=i\circ f\circ j:(S^2,P)\to (S^2,P)$ is a Thurston map with distinguished set $P$, and this Thurston map has hyperbolic orbifold. We will prove that the virtual endomorphism $\phi_F:G_P\dashrightarrow G_P$ associated to $F$  is not injective. It will then follow from Lemma \ref{prepost} that the virtual homomorphism $\phi_f:G_B\dashrightarrow G_A$ is not injective, concluding the proof. 

The proof of this theorem will use both analytic and algebraic arguments. 
From this, we will derive 

\begin{thm}
\label{thm:dense}
Let $f:(S^2,A)\to (S^2,B)$ be an admissible cover with $\#A=\#B$, and suppose that $f$ has potentially hyperbolic orbifold. Then the closure in the Thurston compactification of each nonempty fiber of $\sigma_f:\TTT_B\to\TTT_A$ map contains the Thurston boundary of $\TTT_B$. 
\end{thm} 

\subsection{Proof of Theorem \ref{thm:noninjective}} 

We begin with some preliminary results of independent interest. 

Suppose $f:(S^2,A)\to (S^2,B)$ is an admissible covering map with virtual homomorphism $\phi_f: G_B \dashrightarrow G_A$; as before, let $H_f<G_B$ be the domain of $\phi_f$. Let $d_B$ be the Teichm\"uller distance on $\TTT_B$, and let $d_A$ be the Teichm\"uller distance on $\TTT_A$. 

Recall that every element $g$ of $G_B$ has a {\em minimal displacement} given by 
\[ \delta_B(g):= \inf_{\tau \in \TTT_B} d_B(\tau, g\cdot\tau) \in [0,\infty).\]
Define the analogous object $\delta_A$ for $G_A$. 
Also recall that we have the functional identity 
\[ \sigma(h\cdot\tau) = \phi_f(h)\cdot\sigma(\tau),\quad \forall\; h \in H_f.\]

The following results are related; however, Proposition \ref{hyplength1} and Corollary \ref{cor:twisttotwist} involve nondynamical statements, while Corollary \ref{hyplength2} and Corollary \ref{cor:notconj} involve dynamical statements. 
\begin{prop}\label{hyplength1}
Let $f:(S^2,A)\to (S^2,B)$ be an admissible cover. For each $h \in H_f$ we have 
\[ \delta_A(\phi_f(h)) \leq \delta_B(h).\]
\end{prop}
\pf Fix $\epsilon>0$.  We have 
\renewcommand{\arraystretch}{2}
\[
\begin{array}{cclr}
\delta_B(h) & = & \inf_\xi d_B(\xi, h\cdot\xi) & \mbox{definition}\\
\; & > & d_B(\tau, h\cdot\tau) - \epsilon & \mbox{for some $\tau\in\TTT_B$} \\
\; & \geq & d_A(\sigma_f(\tau), \sigma_f(h\cdot\tau)))-\epsilon & \mbox{$\sigma_f$ is distance nonincreasing}\\ 
\; & = & d_A(\sigma_f(\tau), \phi_f(h)\cdot\sigma(\tau))-\epsilon  & \mbox{functional identity} \\
\; & \geq & \inf_\xi d_A(\xi, \phi_f(h)\cdot\xi)-\epsilon & \mbox{definition}\\ 
\; & = & \delta_A(\phi_f(h))-\epsilon & \mbox{definition}.
\end{array}
\]
\qed
\begin{cor}\label{hyplength2} Let $f:(S^2,P)\to (S^2,P)$ be a Thurston map with hyperbolic orbifold.
\begin{enumerate}
\item There exists $k\in\N$ such that if $h \in H_{f^{\circ k}}$ is pseudo-Anosov, then 
\[ \delta_P(\phi^{\circ k}_f(h)) < \delta_P(h).\] 
\item If $\sigma_f:\TTT_P\to\TTT_P$ is strictly distance-decreasing and $h \in H_f$, then 
\[ \delta_P(\phi_f(h)) < \delta_P(h).\]
\end{enumerate}
\end{cor}
\pf Being pseudo-Anosov implies we can take $\tau$ to realize the infimum in the definition of $\delta_P(h)$, so that $\epsilon$ in the previous proof is equal to $0$. By \cite[Cor. 3.4]{DH1}, there exists a $k\in\N$ such that the $k$th iterate $\sigma_f^{\circ k}$ strictly decreases Teichm\"uller distances.  
\qed

\begin{cor}
\label{cor:twisttotwist}
Let $f:(S^2,A)\to (S^2,B)$ be an admissible cover. Then the virtual homomorphism $\phi_f:G_B\dashrightarrow G_A$ sends multitwists to multitwists.
\end{cor}

\pf There are no elliptic pure mapping class elements, so twists $h$ are characterized by the condition that $\delta_B(h)=0$ (or $\delta_A(h)=0$ if $h\in G_A$).  But $\delta_B(h) \geq \delta_A(\phi_f(h))$ so  $\delta_B(h)=0 \implies \delta_A(\phi_f(h)) = 0$.   
\qed

\begin{cor}
\label{cor:notconj}
Let $f:(S^2,P)\to (S^2,P)$ be a Thurston map, and $k$ as in Corollary \ref{hyplength2}. If the orbifold of $f$ is hyperbolic and $h \in H_{f^{\circ k}}$ is a pseudo-Anosov element, then the mapping class $\phi^{\circ k}_f(h)$ cannot be conjugate to $h$ in $G_P$.  In particular $\phi^{\circ k}_f(h) \neq h$ and $\phi_f(h) \neq h$.  If $\sigma_f:\TTT_P\to\TTT_P$ is distance decreasing, then the mapping class $\phi_f(h)$ cannot be conjugate to $h$; in particular $\phi_f(h) \neq h$. 
\end{cor}

\noindent{\bf Remark.}  One can show Proposition \ref{hyplength1} analytically in another way, via the correspondence on moduli space and the definition of the virtual endomorphism in terms of the induced map on fundamental group determined by the correspondence.  As before, choose basepoints $w_\circledast, \mu_\circledast, \tilde{\mu}_\circledast$ with $Y(w_\circledast)=\mu_\circledast, X(w_\circledast)=\tilde{\mu}_\circledast$. Represent the mapping class $h\in G_B$ by a loop $\gamma$ in $\MMM_B$ based at $\mu_\circledast$.  By suitable choice of basepoints we may assume the length of the loop $\gamma$ is very close to $\delta_B(h)$.  If $h\in H_f$, this loop lifts to a loop $\tilde{\gamma}$ based at $w_\circledast$; so the length of $\tilde{\gamma}$ and of $\gamma$ are the same.  The map $X:\WWW_f\to \MMM_A$ cannot increase lengths, so we conclude that $X(\tilde{\gamma})$ in $\MMM_A$ under this projection has length which is at most that of $\gamma$ in $\MMM_B$.  
\gap

\pf (of Theorem \ref{thm:noninjective}) 
Let $S_m$ denote the $m$-times-punctured sphere. We will need a few results about its extended (nonpure) mapping class group, $\Mod(S_m)$.  
The first is the following result of Bell and Margalit \cite[Thm. 1]{bell:margalit:cohopfian}.   
\begin{thm}
\label{thm:cohopfian} 
Let $m \geq 5$. If $H$ is a finite index subgroup of $\Mod(S_m)$, and $\rho: H \to \Mod(S_m)$ is an injective homomorphism, then there is a unique $g \in \Mod(S_m)$ so that $\rho(h)=ghg^{-1}$ for all $h \in H$.
\end{thm}
Recall we are proving that if $f:(S^2,P)\to (S^2,P)$ is a Thurston map with hyperbolic orbifold, then the kernel of the virtual endomorphism $\phi_f:G_P\dashrightarrow G_P$ contains a pseudo-Anosov element. 
Note that $\Mod(S^2, P)=G_P$ is canonically identified with a finite index subgroup of $\Mod(S_m)$.

Assume first that $\#P \geq 5$.  

Suppose to the contrary that  $\phi_f$ is injective.   Let $k$ be the integer of Corollary \ref{hyplength2}. Then $\phi_f^{\circ k}$ is also injective. By Theorem \ref{thm:cohopfian} applied to $\rho=\phi_f^{\circ k}$ and $H:=H_{f^{\circ k}}$, there is a unique $g \in G_P$ so that $\phi^{\circ k}_f(h)=ghg^{-1}$ for all $h \in H$.  Since $H$ is a finite-index subgroup of $G_P$, it contains a pseudo-Anosov element, $h$.  But then $\phi^{\circ k}_f(h)=ghg^{-1}$, contradicting Corollary \ref{cor:notconj}. 
Let $N:=\ker(\phi_f)$ and let $H$ now denote the domain of $\phi_f$.  We have just shown that $N$ is nontrivial. 


Now assume that $\#P=4$, so that $\MMM_P$ and $\WWW_f$ are each Riemann surfaces of finite type. Assuming that  $\phi_f$ is injective, Theorem \ref{thm:bdry} implies $X:\WWW_f\to\MMM_P$ is proper and surjective; it is also open. It follows that $X: \WWW_f \to \MMM_P$ is a finite branched covering. Let $w_\circledast\in\WWW_f$, and $m_\circledast:=X(w_\circledast)$ be basepoints, and consider $X_*:\pi_1(\WWW_f,w_\circledast)\to\pi_1(\MMM_P,m_\circledast)$. Suppose $X_*$ is injective. 

Let $p: \ZZZ \to \MMM_P$ be the covering of $\MMM_P$ induced by the image subgroup $X_*(\pi_1(\WWW_f,w_\circledast))$.  By elementary covering space theory, there exists a lift $\tilde{X}: \WWW_f\to \ZZZ$ of $X$ satisfying $p\circ\tilde{X}=X$.   Since degrees of branched coverings are multiplicative, $\tilde{X}$ is a finite branched covering. By construction, the induced map on fundamental groups $\tilde{X}_*$ is injective and surjective, hence an isomorphism.  Since $\MMM_P$ is a thrice-punctured sphere, $\ZZZ$ and $\WWW_f$ are punctured surfaces of strictly negative Euler characteristic, say $-k_\ZZZ, -k_{\WWW_f}$, respectively. The fundamental groups of $\ZZZ$ and $\WWW_f$ are then free groups of ranks respectively $1+k_\ZZZ, 1+k_{\WWW_f}$. Since they are isomorphic (via $X_*$), these ranks coincide, hence so do the Euler characteristics of $\ZZZ$ and $\WWW_f$. By the Riemann-Hurwitz formula, the degree of $\tilde{X}$ must be $1$. It follows that $\tilde{X}$ is an unramified cover, and hence that $X=p\circ\tilde{X}$ is an unramified cover as well.\footnote{The authors acknowledge Allan Edmonds for providing the above arguments in this paragraph.}  But then local inverse branches of $X\circ Y^{-1}$ are isometries with respect to the hyperbolic (equivalently, Teichm\"uller) metric.  Taking a composition of two such branches, we obtain again an isometry. This is impossible, since the orbifold of $F$ is assumed hyperbolic, so that $\sigma_f$ and, hence, that  $X\circ Y^{-1}$ have second iterates that are contractions.

We conclude that in both cases $N$ is nontrivial. Since $S_m$ is torsion-free,  $N$ is infinite. The next results we need concern the notion of irreducibility: a subgroup $H < \Mod(S_m)$ is {\em reducible} if there exists a nonempty multicurve which is stabilized by every element of $H$.   A subgroup containing a pseudo-Anosov element cannot be reducible.  Hence any finite-index subgroup $H< \Mod(S_m)$ is irreducible, since given any pseudo-Anosov element, some power will be a pseudo-Anosov element which lies in $H$.   By \cite[Cor. 7.13]{ivanov:book:subgroups}, the subgroup $N$ is irreducible, and by   \cite[Cor. 7.14]{ivanov:book:subgroups}, $N$  contains a pseudo-Anosov element.\footnote{The authors acknowledge Dan Margalit for this reference.}
\qed 

\subsection{Proof of density, Theorem \ref{thm:dense}} 
We will prove that if $f:(S^2,A)\to (S^2,B)$ is an admissible map with $\#A=\#B$ and potentially hyperbolic orbifold, then the closure in the Thurston compactification of each nonempty fiber of $\sigma_f:\TTT_B\to\TTT_A$ contains the Thurston boundary of $\TTT_B$. 

Let $N=\ker(\phi_f)$; by Theorem \ref{thm:noninjective}, $N$ contains pseudo-Anosov elements.   Suppose $\tau \in \TTT_P$ is an element of a nonempty fiber $E$ of $\sigma_f$.    Given a subgroup of $G_B$, recall that its {\em limit set} in the Thurston boundary of $\TTT_B$ is defined as the closure of the set of fixed-points of its pseudo-Anosov elements. 
For $n \in N$ we have $\sigma_f(n\cdot\tau) = \sigma_f(\tau)$, so that $E$ is $N$-invariant. Since pseudo-Anosov elements have north-south dynamics on the Thurston compactification of ${\TTT}_P$, the accumulation set of $E$ in the Thurston compactification contains the limit set $\Lambda(N)$. But $N$ is a normal subgroup of $H_f$ so if $x\in \Lambda(N)$ is a fixed-point of some $n\in N$, then $h\cdot x$ is a fixed-point of $h\circ n \circ h^{-1} \in N$ for all $h \in H_f$. Thus $\Lambda(N)$ is $H_f$-invariant and it follows that $\Lambda(N) = \Lambda(H_f)$. Since $H_f$ is of finite index in $G_B$, we conclude  $\Lambda(N) = \Lambda(H_f)=\Lambda(G_B)$, which is equal to the Thurston boundary of $\TTT_B$. 
\qed

\section{Finite global attractor for the pullback relation} 
\label{section:FGA}
For the remainder of the article, we restrict to the dynamical setting; that is, $f:(S^2,P)\to (S^2,P)$ is a Thurston map with postcritical set $P$. In this section, we establish conditions on $f$ under which there exists a finite set $\AAA \subset \SSS_P$  such that $\AAA \pullback \AAA$ and the orbit of each curve eventually falls into $\AAA$; we call $\AAA$ a {\em finite global attractor}. 

In this paragraph, we state a known algebraic condition for the existence of a finite global attractor. Suppose $S$ is a finite generating set for a group $G$.  We denote by $||g||$
the word length of $g$ in the generators $S$.
\begin{defn}
\label{defn:contracting1}
The virtual endomorphism $\phi: G  \dashrightarrow G$ is called {\em contracting}
if the {\em contraction ratio}
\[ \rho:= \limsup_{n\to\infty}\left( \limsup_{||g||\to\infty}
\frac{||\phi^n(g)||}{||g||}\right)^{1/n} < 1.\]
\end{defn}
The contraction ratio of the virtual endomorphism $\phi$ is
independent of the choice of generating set, and is an asymptotic property: $\phi$ is contracting if and only if some iterate of $\phi$ is contracting. 
  In \cite[Thm. 1.4]{kmp:tw}, it is shown that the algebraic criterion of contraction  of the virtual endomorphism $\phi_f:G_P\dashrightarrow G_P$ implies the existence of a finite global attractor for $\SSS_P \pullback \SSS_P$.
It is easy to construct examples of obstructed maps $f$ for which $\SSS_P \pullback \SSS_P$ does not have a finite global attractor---for example, there could be a large essential subsurface $\Sigma \subset S^2\setminus P$ for which $f|_\Sigma = \id_\Sigma$, leading to e.g. infinitely many fixed-curves.  It is perhaps somewhat more surprising that one can achieve this with an expanding map.  
\gap

\noindent{\bf Example:} Let $f$ be the Thurston map obtained by starting with the degree four integral Latt\`es example and applying the ``blowing up'' surgery  along a single vertical arc joining two critical points (see Figure 4).  
The domain and range spheres are each the union of the two squares $A, B$ at the left side of Figure \ref{fig:ex1} along their common boundary. The arrows go in the direction of pulling back. The map $f$ sends small topological quadrilaterals labelled $A, B$ to the large squares labelled $A, B$, respectively; their boundary edges $e_i$ map to the boundary edge with the same label $e_i$.  As usual set $P:=P_f$. 
Then the pullback relation $\SSS_P \pullback \SSS_P$ has infinitely many fixed curves.  To see this, observe that the  horizontal and vertical curves $\gamma_h, \gamma_v$ are fixed.  Since $\lambda_f(\gamma_v)=\gamma_v$, it follows from \cite[Thm. 8.2]{kmp:cds} that if $T^2$ denotes a double right Dehn twist about $\gamma_v$, then $T^2\circ f = f \circ T^2$ up to homotopy relative to $P$. Thus the curves $\gamma_n:=T^{2n}(\gamma_h), n \in \Z$ are also each fixed by $\pullback$.  We remark that since the subdivision rule describing $f$ has mesh tending to zero combinatorially, there exists by \cite[Thm. 2.3]{cfp:fsr} an expanding map $g$ homotopic to $f$ relative to $P$ arising as the subdivision map of the shown subdivision rule. 
\gap

\begin{figure}
\begin{center}
\includegraphics[width=4in, angle=-90]{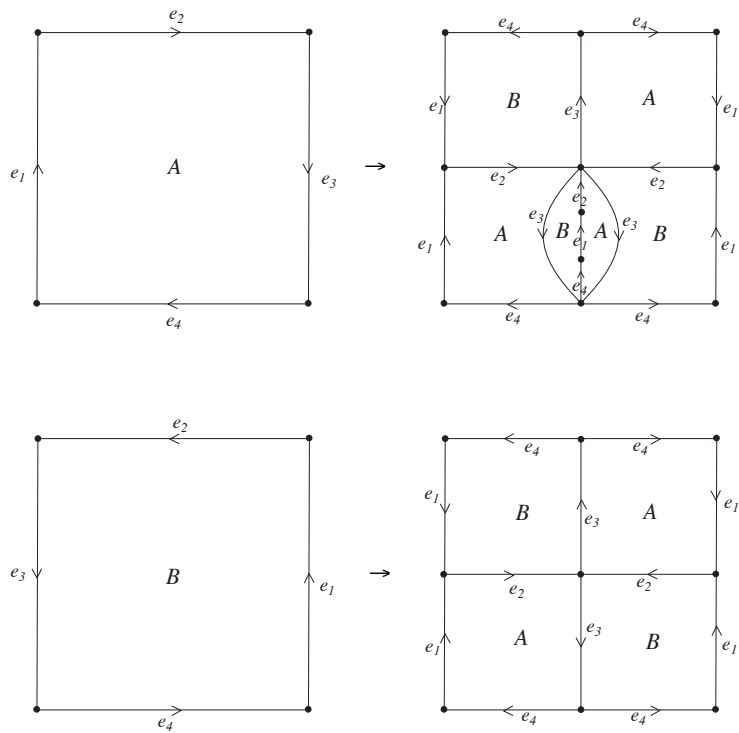}
\caption{(Image by W. Floyd.) The horizontal curve is fixed under pullback, as are all elements of its orbit under iteration of a double Dehn twist about the vertical curve.}
\end{center}
\end{figure}

It is therefore natural to restrict attention to the case when $f$ is a rational map, in which case $\sigma_f:\TTT_P\to\TTT_P$ has a unique fixed point $\tau_\circledast$; let $m_\circledast=\pi_P(\tau_\circledast)$.  Below, we give two different analytic conditions on the correspondence $X\circ Y^{-1}: \MMM_P \to \MMM_P$  which imply that $\SSS_P \pullback \SSS_P$ has a finite global attractor. These properties depend only on the $(A,B)$-Hurwitz class of $f$ and not on the choice of base fixed-point. 

{In preparation for stating our results, the following two paragraphs give definitions. }

{A nonempty subset $K \subset \MMM_P$ is {\em invariant} under the correspondence $X \circ Y^{-1}$ if $X\circ Y^{-1}(K) \subset K$.  Consider the case when there exists a nonempty compact invariant subset, $K$.  This condition is quite strong.   If $K$ is such a set, then so is any $R$-neighborhood of $K$, defined using the Teichm\"uller metric, since lengths of paths do not increase under application $X \circ Y^{-1}$.  In this case it  follows that $\MMM_P$ is exhausted by compact invariant connected sets.  Informally, then, the condition of having a nonempty invariant compact subset may be thought of as asserting that the ends of moduli space are repellors of $X\circ Y^{-1}$.  The relationship between the property of having a nonempty invariant compact subset and other topological, dynamical, and algebraic properties is investigated in \cite{koch:pilgrim:selinger:corlim}.  }

{We will call a length metric $\ell$ on $\MMM_P$  {\em WP-like} }if its lift $\tilde{\ell}$ to $\TTT_P$ has the property that the identity map defines a homeomorphism between the completions of $\TTT$ with respect to $\tilde{\ell}$ and with respect to the WP metric. This implies, in particular, that $(\MMM, \ell)$ has finite diameter. 
We say that the pullback correspondence $X\circ Y^{-1}$ on $\MMM_P$ is {\em uniformly contracting} with respect to an WP-like  length metric $\ell$ on $\MMM_P$ if there is a constant $0 \leq \lambda < 1$ such that for each curve $\gamma: [0,1] \to \MMM_P$ of finite length, and each lift $\tilde{\gamma}$ of $\gamma$ under $Y$, we have  $\ell(X \circ \tilde{\gamma}) \leq \lambda \cdot \ell(\gamma)$, where $\ell(\gamma)$ is the length of $\gamma$.

\begin{thm}
\label{thm:eventually_finite}
Suppose $f:(\P^1,P)\to (\P^1,P)$ is a rational Thurston map with hyperbolic orbifold. 
\be 
\item 
If the correspondence $X\circ Y^{-1}:\MMM_P\rightrightarrows\MMM_P$ has a nonempty compact invariant subset, then $\phi_f$ is contracting, and the pullback relation $\pullback$ has a finite global attractor. 

\item If the correspondence $X\circ Y^{-1}:\MMM_P\rightrightarrows\MMM_P$ is uniformly contracting with respect to a WP-like length metric $\ell$, then under iteration of the pullback relation $\pullback$, every curve becomes trivial, i.e. the set $\{o\}$ consisting of the trivial curves (inessential and peripheral) is a finite global attractor. 
\eb
\end{thm}

\noindent{\bf Examples.}
\be
\item The correspondence induced by the Rabbit polynomial $f(z)=z^2+c_R$ satisfies the hypothesis in $(1)$.  As shown in \cite[\S 5.1]{bartholdi:nekrashevych:twisted}, the map $X$ is an inclusion, and in coordinates identifying $\MMM_P$ with $\IP^1\setminus\{0, 1, \infty\}$, the map $Y\circ X^{-1}$ is given by the restriction of the critically finite hyperbolic rational function $g(w) = 1-\frac{1}{w^2}$ to the complement of $\{0, \infty, \pm 1\}$. The Julia set $K$ of $g$ is then a nonempty compact invariant subset for $X \circ Y^{-1}=g^{-1}$.  Theorem \ref{thm:eventually_finite}(1) implies that the pullback relation on curves has a finite global attractor.  It is computed in \cite[\S 8] {kmp:tw} using algebraic techniques. 

The same arguments apply in the case of preperiod $2$, period $1$ quadratic polynomials; see \cite[\S 7.2]{bartholdi:nekrashevych:twisted} and \cite[\S 9]{kmp:tw}.  

\item The correspondence induced by the dendrite Julia set polynomial $f(z)=z^2+i$ satisfies the hypothesis in $(2)$. As shown in \cite[\S 6.2]{bartholdi:nekrashevych:twisted}, the map $X$ is an inclusion, and in coordinates identifying $\MMM_P$ with $\IP^1\setminus\{0, 1, \infty\}$, the map $Y\circ X^{-1}$ is given by the restriction of the critically finite Latt\`es-type subhyperbolic (but not hyperbolic) rational function $g(w) = (1-2/w)^2$ to the complement of $\{0, \infty, 1, 2\}$.   The map $g$ defines a self-orbifold-cover of the Euclidean $(2,4,4)$-orbifold (see \cite[\S 9]{DH1}) which uniformly expands the corresponding Euclidean metric $\ell$ by the factor $|1+i|=\sqrt{2}$.  The identity map of $\MMM_P$ defines a homeomorphism between the completions with respect to the Weil-Petersson and Euclidean metric, and it follows easily that $\ell$ is WP-like. Theorem \ref{thm:eventually_finite}(2) implies that under iterated pullback, each curve becomes trivial. This is established in \cite[\S 9]{kmp:tw} using algebraic techniques. Guizhen Cui (personal communication) has pointed out that there is an elementary proof of this fact: it is easy to see that any curve which does not surround the unique 2-cycle in $P_f$ must become trivial, and any curve which does must, by the Schwarz-Pick lemma, have geodesic representative which strictly shrinks under pullback.  

\item The correspondence induced by $f(z)=\frac{3z^2}{1z^3+1}$, studied in \cite[\S 4]{bekp} satisfies the hypotheses of neither (1) nor (2) above.  The analysis in \cite{lodge:thesis}, however, shows that nonetheless the pullback relation on curves has a finite global attractor. 
\eb

\pf Since $f$ is assumed rational, $\sigma_f:\TTT_P\to\TTT_P$ has a fixed-point which we denote by $\tau_\circledast$; it projects to $w_\circledast\in\WWW_f$ and $m_\circledast\in\MMM_P$.
  
We first prove $(1)$ by showing that the virtual endomorphism $\phi_f:G_P\dashrightarrow G_P$  is contracting, and appealing to \cite[Theorem 1.4]{kmp:tw}.  To see that $\phi_f$ is contracting, recall that $\phi_f$ is also the virtual endomorphism on $\pi_1(\MMM_P,m_\circledast)$ given by $(X\circ Y^{-1})_*=X_*\circ Y^*$, where  $Y^*(\gamma)$ is the lift under $Y$ of a loop $\gamma \in \pi_1(\MMM_P,m_\circledast)$ based at $w_\circledast\in\WWW_f$.  Equip $\MMM_P$ with the Teichm\"uller metric (equivalently, by Royden's theorem, the Kobayashi metric) and $\WWW_f$ with the lift of this metric.   
Given a smooth loop in $\MMM_P$ based at $m_\circledast$ in the domain of the virtual endomorphism, it lifts under $Y$ to a loop of the same length in $\WWW_f$ and projects to a loop in $\MMM_P$, whose length is strictly shorter.   The space $\MMM_P$ is a hyperplane complement, so its fundamental group is generated by a finite collection $S$ of loops.  The union of these loops lies in some invariant compact subset $K$; we may assume $K$ is a compact submanifold with smooth boundary.  In $K$, the word length $||g||_S$ with respect to $S$ is comparable to the infimum of lengths of loops in $K$ representing $g$. 
Since $K$ is compact and invariant, and since the orbifold of $f$ is hyperbolic,  the second iterate of  the correspondence is uniformly contracting on $K$.  It follows that lengths of such liftable loops in $K$, and hence the word lengths of the corresponding group elements,  are uniformly contracted under iteration and so $\phi_f$ is contracting. 

We now prove $(2)$.  Let $\tilde{\ell}$ denote the lifted WP-like metric. Results of \cite[\S 4]{selinger:boundary} imply that the pullback map $\sigma_f$ extends to the Weil-Petersson completion $\overline{\TTT}_P$, which, by the WP-like property of $\ell$, coincides with the $\tilde{\ell}$ completion of $\TTT_P$.   The hypothesis implies that $\sigma_f$ is uniformly contracting with respect to $\tilde{\ell}$, and hence that under iteration, each point in $\overline{\TTT}_P$ converges exponentially fast to $\tau_\circledast$ under iteration. 
Now let $\Gamma$ be a multicurve on $S^2\setminus P$. Proposition \cite[Prop. 6.1]{selinger:boundary} implies that for each multicurve $\Gamma$, the pullback map $\sigma_f$ sends the stratum $\TTT^\Gamma_P$ to the stratum $\TTT^{f^{-1}(\Gamma)}_P$.    

Define the {\em complexity} of $\Gamma$, denoted $\kappa(\Gamma)$, to be the distance between $\tau_\circledast$ and the stratum $\TTT^\Gamma_P$ with respect to the lifted metric $\tilde{\ell}$. 
The admissibility of $\ell$ implies 
\be
\item for all $\Gamma$, $\kappa(\Gamma) < \infty$, and 
\item there exists $r_0>0$ such that for all $\Gamma \neq \emptyset$, $\kappa(\Gamma)>r_0$.
\eb
The hypothesis, the admissibility of $\ell$, and the fact that $\sigma_f(\TTT^\Gamma_P) \subset \TTT^{f^{-1}(\Gamma)}_P$ implies 
\be
\setcounter{enumi}{2}
\item there exists $\lambda < 1$ such that for all $\Gamma$, $\kappa(f^{-1}(\Gamma)) \leq \lambda\cdot\kappa(\Gamma)$. 
\eb
From $(1), (2)$, and $(3)$, it follows that given any multicurve $\Gamma$, there exists $n \in \N$ such that $\kappa(f^{-n}(\Gamma))=0$ and hence that $f^{-n}(\Gamma)$ is trivial.   

Hence under iterated pullback, each multicurve becomes trivial. It follows that every curve becomes trivial under iterated pullback. 

\qed 

\section{Shadowing} 
In this section, we briefly consider one consequence of surjectivity of the virtual endomorphism; this will be used in \S 9.2 below. 

\begin{prop}
\label{prop:shadowing}
Let $f:(S^2,P)\to (S^2,P)$ be a Thurston map. Suppose that $\phi_f:G_P\dashrightarrow G_P$ is surjective. Then any finite orbit segment of the correspondence $X\circ Y^{-1}:\MMM_P\rightrightarrows\MMM_P$ is covered by an orbit segment of the pullback map $\sigma_f:\TTT_P\to\TTT_P$.
\end{prop}

\pf Let $\mu_1\in\MMM_P$, and let $\mu_2\in X(Y^{-1}(\mu_1))\subset\MMM_P$. By the Fundamental Diagram, Figure \ref{fund}, there exist $\tau_1',\tau_2'\in \TTT_P$ so that $\sigma_f(\tau_1')=\tau_2'$, and $\pi_P(\tau_1')=\mu_1$ and $\pi_P(\tau_2')=\mu_2$. Let $\tau_2\in\pi_P^{-1}(\mu_2)$. Then there exists $g\in G_P$ so that $\tau_2=g\cdot\tau_2'$. Since $\phi_f:G_P\dashrightarrow G_P$ is surjective, there exists $h\in G_P$ so that $\phi_f(h)=g$. Define $\tau_1:=h\cdot \tau_1$. Then 
\[
\sigma_f(\tau_1) = \sigma_f(h\cdot\tau_1') = \phi_f(h)\cdot\sigma_f(\tau_1') = g\cdot\tau_2' = \tau_2.
\]
The claim now readily follows by induction.

\section{Obstructed twists and repelling fixed points in $\bdry \MMM_P$} 

If $f$ is a Thurston  map with $\#P=4$, then the projection maps $Y, X: \WWW_f \to \MMM_P$ extend to holomorphic maps $Y, X: \overline{\WWW}_f \to \overline{\MMM}_P \simeq \IP^1$ yielding a holomorphic correspondence on the compactified moduli space $X\circ Y^{-1}: \IP^1 \to \IP^1$.

\subsection{Obstructed twists} 
\begin{thm} 
\label{thm:tfae} 
Let $f:(\P^1,P)\to (\P^1,P)$ be a rational Thurston map  with hyperbolic orbifold. Suppose $\#P=4$. 

Then the following conditions are equivalent.
\be
\item There exists $g \in G_P$, a twist $T$, and a nonzero $k\in\Z$ such that $f_*:=g \circ f$ commutes up to homotopy with $T^k$.
\item There exists $g \in G_P$ such that $f_*:=g\circ f$ has an obstruction $\Gamma:=\{\gamma\}$ for which $\lambda_f(\gamma)=\gamma$.
\item There exists an element $g \in G_P$, a twist $T$, and a nonzero integer $k$ such that $(T^k)^g\in H_f$ and $\phi_f((T^k)^g)=T^k$.
\item There exist curves $\gamma_1, \gamma_2$ in $\P^1\setminus P$ such that (i) $g(\gamma_1) =\gamma_2$ for some $g \in G_P$, and (ii) $\lambda_f(\gamma_1) = \gamma_2$.
\item There exists a point $0 \in \bdry \MMM_P$ and a single-valued branch $\mu \mapsto \beta(\mu)$ of the pullback correspondence $X \circ Y^{-1}$ near $0$ such that $\beta(0)=0$ and, in suitable local holomorphic coordinates for which $0$ is the origin, $\beta(\mu) = a\mu + O(\mu^2)$, where $0 < |a| < 1$. 
\eb
\end{thm}

\pf  $(1) \iff  (2)$  We take the elements $g$ to be the same. 

Condition (1) is equivalent to the condition that $T^k \circ f_* = f_* \circ T^k$ up to homotopy, i.e. that $T^k \in \dom(\phi_{f_*})$ and $\phi_{f_*}(T^k)=T^k$. 
In turn this is equivalent to the condition that for some $k'$ with $|k'| \geq |k|$ we have $T^k \in \dom(\phi_f')$ and $\phi_{f_*}(T^{k'})=T^{k'}$.  Denote by $\gamma$ the core curve of $T$.   Put $w=k'\gamma$ so that the multitwist $M_w:=T^{k'}$. 

If Condition (1) holds, then equation (\ref{eqn:philambda}) implies  
\[ \phi_{f_*}(M_w) =M_w \implies \lambda_{f_*}(w)=w \implies \lambda_{f_*}(k'\gamma)=k'\gamma \implies \lambda_{f_*}(\gamma)=\gamma\]
by linearity, and so $\gamma$ is an obstruction for $f_*$  with eigenvalue $1$. 

Now suppose condition (2) holds, with the obstruction being given by $\gamma$. Choose $k' \in \Z$ so $T^{k'} \in \dom(\phi'_{f_*})$. Again by equation (\ref{eqn:philambda}) we have 
\[ \lambda_{f_*}(\gamma)=\gamma \implies \lambda_f(k'\gamma)=k'\gamma \implies \phi_{f_*}(T^{k'})=T^{k'}\]
as required. 

$(1) \iff (3)$  By applying Lemma \ref{prepost} with $i=g$ and $j=\id$, for any $g \in G$ we have 
\[ \phi_f((T^k)^g) = T^k \iff \phi_{g \circ f}(T^k) = T^k.\]

$(4) \iff (2)$  This follows immediately from the observation that $\lambda_{g\circ f}(\gamma) = \lambda_f(g^{-1}(\gamma))$.  

$(1) \implies (5)$  If $T^k \circ f_* = f_* \circ T^k$ up to homotopy then $\sigma_{T^k}\circ \sigma_{f_*} = \sigma_{f_*} \circ \sigma_{T^k}$.  Let $\gamma$ be the core curve of the twist $T$, let $\Gamma = \{\gamma\}$, and consider the quotient
\[ \IH/\langle z \mapsto z+2k\rangle \simeq {\TTT}_P/k\Tw(\Gamma)\simeq \D^*\] 
which we identify with the punctured disk. Then $\sigma_f$ descends to an analytic map 
\[ \overline{\sigma}_{f_*}: \D^* \to \D^*.\]
The fundamental diagram becomes 
\[\xymatrix{  & \TTT_P \ar[d]\ar[rr]^{\sigma_{f_*}} & &
\TTT_P \ar[d] \\
& \D^* \ar[dd]_{\overline{\pi}_P}\ar[rr]^{\overline{\sigma}_{f_*}} \ar[dr]^{\overline{\omega}_f}& &
\D^*\ar[dd]^{\overline{\pi}_P} \\
&&\mathcal{W}_f\ar[dl]_Y \ar[dr]^X &\\
& \MMM_P & & \MMM_P.}
\]
where the vertical arrows are covering maps. 
Since $\Gamma$ is $f_*$-invariant,  the induced map $\overline{\sigma}_{f_*}$ extends over the origin, yielding an analytic map 
\[  \overline{\sigma}_{f_*}: (\D,0) \to (\D, 0).\]
For convenience of notation, we denote by $0$ the corresponding end of $\MMM_P$ and of $\WWW_f$.  We obtain a commutative diagram of analytic maps 
\[\xymatrix{  
& (\D, 0) \ar[dd]_{\pi_P}\ar[rr]^{\overline{\sigma}_{f_*}} \ar[dr]^{\overline{\omega}_f}& &
(\D, 0)\ar[dd]^{\pi_P} \\
&&(\mathcal{W}_f\union\{0\}, 0)\ar[dl]_{Y} \ar[dr]^{X} &\\
& (\MMM_P\union\{0\}, 0) & & (\MMM_P\union\{0\}, 0).}
\]
Since $f$ has hyperbolic orbifold, $\overline{\sigma}_{f_*} \not\in\Aut(\D)$.  By the Schwarz Lemma, its derivative at the origin therefore satisfies $|(\overline{\sigma}_{f_*})'(0)| < 1$.  Taking inverse branches near the origin, we have 
\[ (X\circ Y^{-1})'(0) = (\pi_P \circ \overline{\sigma}_{f_*} \circ \pi_P^{-1})'(0)=\sigma_{f_*}'(0) \]
and so the pullback correspondence $X \circ Y^{-1}$ has an attracting fixed point in the boundary of moduli space. 

In the remaining  paragraphs, we prove that $\sigma_{f_*}'(0) \neq 0$, and that $(5) \implies (3)$. 

We first discuss some specializations of the results in the ``Fivefold Way'' Theorem \ref{thm:fivefold_way} to the case $\#P=4$; cf. \cite{lodge:thesis}. We may identify $P=P_f=\{0, 1, \infty, m_\circledast\}$, where $m_\circledast\in X\circ Y^{-1}(m_\circledast)$ (since $f$ is rational). Recall that there is a natural identification $G_P \to \pi_1(\MMM_P, m_\circledast)$ such that $\phi_f: G_P\dashrightarrow G_P$  coincides with the virtual endomorphism $X_*\circ Y^* $ of $\pi_1(\MMM_P, m_\circledast)$ induced by the correspondence $X\circ Y^{-1}: (\MMM_P, m_\circledast)\rightrightarrows (\MMM_P, m_\circledast)$.   Under this identification, a simple oriented loop $\gamma_{loop}$ based at $m_\circledast$ in $\MMM_P$ which is peripheral and which encloses $0$ on its left-hand side corresponds to the right Dehn twist $T$ about the right-hand boundary component $\gamma_{curve}$ of a regular neighborhood of $\gamma_{loop}$ in $\IP^1\setminus \{0, 1, \infty, m_\circledast\}$. Denote by $\lambda$  the  group element of $\pi_1(\MMM, m_\circledast)$ corresponding to $g$. We see that
\[ \phi_{g \circ f}(T^k)=T^k \iff \phi_f((T^k)^g)=(T^k)^g \iff (X_*\circ Y^*)(\lambda\cdot\gamma^k_{loop}\cdot\lambda^{-1})=\gamma^k_{loop}  \]
where $\cdot$ denotes concatenation of paths, $\lambda^{-1}$ is the path $\lambda$ traversed in the opposite direction, and in the right-hand expression the path $\lambda$ is traversed first. 

Now suppose $(1)$ holds, so that $(X_*\circ Y^*)(\lambda\cdot\gamma^k_{loop}\cdot\lambda^{-1})=\gamma^k_{loop} $. 
The assumption that $\lambda\cdot\gamma^k_{loop}\cdot\lambda^{-1}$ is in the domain of $(Y^*)$ implies that the $k$-th power of a tiny peripheral simple closed curve surrounding $0\in\MMM_P$ lifts under $\overline{Y}$ to a tiny simple closed curve surrounding $0 \in \WWW_f$.  Since $(X_*\circ Y^*)(\lambda\cdot\gamma^k_{loop}\cdot\lambda^{-1})=\gamma^k_{loop} $, the tiny curve surrounding $0\in\WWW_f$, obtained by lifting, maps by $X$ to the $k$-th power of a simple closed curve surrounding $0\in\MMM_P$. This shows that the local degree of the branch of $X\circ Y^{-1}$ at $(0,0)$ is 1 and, hence, the derivative at this point does not vanish.
This finishes the proof of $(1) \implies (5)$.

The proof of $(5) \implies (3)$ proceeds similarly.  If $(5)$ holds, then the local degrees of $\overline{X}$ and $\overline{Y}$ at $0$ are the same, say $k$. Thus under the correspondence a $k$-th power of a tiny simple peripheral curve about $0$ lifts under $\overline{Y}$ and projects under $\overline{X}$ to a $k$-th power of a tiny simple peripheral curve about $0$. Since free homotopy classes of curves correspond to conjugacy classes of loops based at $m_\circledast$, we conclude that $(X_*\circ Y^*)(\lambda\cdot\gamma^k_{loop}\cdot\lambda^{-1})=\gamma^k_{loop}$ for some loop $\gamma_{loop}$ surrounding $0$ and some loop $\lambda$ based at $m_\circledast$. Thus $(3)$ holds, and the proof is complete. 
\qed

\subsection{Dynamical consequences} 

\begin{thm}
\label{thm:dyncons} 
Under the hypotheses and notation of Theorem \ref{thm:tfae}, assume further that the map $X$ of the correspondence on moduli space is injective, so that an inverse of the pullback map descends and extends to a holomorphic self-map\footnote{A so-called ``$g$-map'', in the language of the first author.} $Y\circ X^{-1}:  \IP^1 \to \IP^1$.  Then the following further properties hold. 
\be[I.]
\item For any Thurston map $F:(S^2,P)\to (S^2,P)$ which is $(P,P)$-Hurwitz equivalent to $f$, the associated virtual endomorphism $\phi_F$ on $G_P$ is surjective. 
\item The end $0$ is a repelling fixed-point of $Y \circ X^{-1}$. 
\item Conditions $(1)$ and $(3)$ hold with $k=1$. 
\item The univalent pullback relation on curves induced by the rational map $f$ has arbitrarily long nontrivial orbits (necessarily comprised of distinct curves), whose elements are all homologous (and hence differ by elements of $G_P$) .
\item The collection of obstructed Thurston maps $T^n\circ f_*$, $n \in \Z$, are pairwise combinatorially inequivalent.  
\item For any $\delta>0$ and any sequence $\epsilon_1>\epsilon_2> \ldots \epsilon_k$ there exists $\tau_0 \in \TTT_P$ and iterates $n_1>m_1>n_2>m_2 > \ldots > n_k>m_k$ such that (i) $\pi_P(\tau_0)=m_\circledast$; (ii) for each $i=1, \ldots, k$,
\[ \ell(\tau_{n_i}) < \epsilon_i, \;\; d(\pi_P(\tau_{m_i}), m_\circledast) < \delta\]
where $\ell(\tau)$ is the length of the systole on the underlying Riemann surface. 
\eb
\end{thm} 

\pf 

\be[I.]
\item Different choices of representatives in the $(P,P)$-Hurwitz class of $f$ lead to virtual endomorphisms which differ by inner automorphisms.  So (just to fix ideas) it suffices to verify this for the case of the rational map, $f$, corresponding to the fixed-point $m_\circledast$.  But this is clear, since the nonconstant inclusion $X: \WWW_f \hookrightarrow \MMM_P$ of finitely punctured spheres must induce a surjection on fundamental groups. 

\item  This follows immediately from condition $(5)$ in Theorem \ref{thm:tfae}. 

\item The previous conclusion implies that $\overline{Y}$ is locally injective near $0$, and hence that a simple peripheral loop $\gamma_{loop}$ in $\MMM_P$ based at $m_\circledast$ surrounding the end $0$ in $\MMM_P$ on its left-hand side lifts under $Y$ to a simple loop in $\WWW_f$ based at $w_\circledast$ surrounding the end $0$ in $\WWW_f$. This implies that the corresponding twist $T$ lies in $H_f$. Since $X$ is injective and the end $0$ in $\WWW_f$ maps to the end $0$ in $\MMM_P$, we have that $\phi_f(T)$ is conjugate to $T$ via some $g \in G_P$, i.e. $\phi_f(T^g)=T$ and so condition $(3)$ of Theorem \ref{thm:tfae} holds with $k=1$. The proof of the equivalence of $(1)$ and $(3)$ in Theorem \ref{thm:tfae} shows that the integers $k$ in $(1)$ and in $(3)$ are the same. 

\item The argument given below is due to L. Bartholdi.  By the previous assertion, we may suppose condition $(3)$ in Theorem \ref{thm:tfae} holds with $k=1$, and let $f, g, T, \phi_f$ be as in the statement there. Recall that $G_P$ is free, generated by any pair of Dehn twists whose core curves intersect minimally. Choose such a twist $S \in G_P$ so that $G_P=\genby{S,T}$. Since $\phi_f$ is surjective, there exists $\tilde{T}:=T^g$ and $\tilde{S} \in G$ such that $\phi_f(\tilde{T})=T, \phi_f(\tilde{S})=S$.  Since $G_P$ is free on the generators $S, T$, the homomorphism defined by sending $T \mapsto \tilde{T}, S \mapsto \tilde{S}$ is  an injection $\sigma: G_P \to G_P$ giving a section of $\phi_f$, i.e. $\phi_f \circ \sigma = \id_{G_P}$. 
For $n \in \N$ let 
\[w_n = g \cdot \sigma(g) \cdot \ldots \cdot \sigma^{\circ(n-1)}(g).\]
Then an easy induction argument shows 
\[ \phi_f(T^{w_n}) = T^{w_{n-1}}.\]
The elements $T^{w_n}$ are a primitive Dehn twists whose core curve $C_n$ lifts under $f$ to the core curve $C_{n-1}$ of $T^{w_{n-1}}$. These curves must be pairwise distinct (else $f$ is obstructed) and pairwise homologous (since the twists are pairwise conjugate). This proves $(4)$. 


\item We again suppose condition $(3)$ in Theorem \ref{thm:tfae} holds with $k=1$, and let $f, g, T, \phi_f$ be as in the statement there.

\begin{lemma}
\label{lemma:obstructions}
Suppose $F$ is a Thurston map  with hyperbolic orbifold for which $\#P_F=4$. If $F$ is obstructed, it has a unique Thurston obstruction $\Gamma=\{\gamma\}$.  If $h \circ F = F \circ h$ up to homotopy for some $h \in G_P$, then $h(\gamma)=\gamma$, preserving orientation.  
\end{lemma}

Cf. \cite[Prop. 6.10]{bartholdi:nekrashevych:twisted}. 

\pf  By the main result of \cite{kmp:canonical}, there is an obstruction (called the \emph{canonical obstruction} of $F$) that  is disjoint from all other obstructions. But on a sphere with four marked points, any two distinct essential nonperipheral curves intersect.  A homeomorphism representing a pure mapping class that yields an automorphism of $F$ up to homotopy must send an obstruction to an obstruction, hence must fix the unique obstruction, preserving orientation since (by pureness) each disk it bounds must be mapped to itself.
\qed

Suppose $n_1, n_2 \in \Z$ and $T^{n_1}\circ f_*$ is combinatorially equivalent to $T^{n_2}\circ f_*$ via $g \in G_P$, so that up to homotopy relative to $P$, 
\[ g \circ T^{n_1}\circ f_*\circ g^{-1} = T^{n_2}\circ f_*.\]
By construction, both $T^{n_1}\circ f_*$ and $T^{n_2}\circ f_*$ have a common obstruction, $\Gamma=\{\gamma\}$, which is the core curve of the twist, $T$.  By Lemma \ref{lemma:obstructions}, $g(\gamma)=\gamma$, preserving orientation. Hence $g=T^l$ for some $l \in \Z$ and 
\[ g \circ T^{n_1}\circ f_*\circ g^{-1} = T^l \circ T^{n_1}\circ f_* \circ T^{-l} = T^{n_2}\circ f_*.\]
Since $f_*$ commutes with $T$ up to homotopy, we have  
\[ f_* \circ T^{n_1} = f_*\circ T^{n_2}\]
and so $n_1=n_2$ since the right action of $G_P$ on the set $\FFF$ of homotopy classes of Thurston maps with postcritical set $P$ is free (\cite[\S 3]{kmp:tw}, \cite[Prop. 3.1]{kameyama:equivalence}). 

\item Since $Y$ is a covering map, the critical values of the map $\overline{Y}\circ \overline{X}^{-1}$ lie in the set of ends of $\MMM_P$.  Since these ends must map to themselves, we conclude that $Y\circ X^{-1}$  is analytically conjugate to a subhyperbolic rational map $R$ with a repelling  fixed point at the origin, $0$; the point $m_\circledast$ is also fixed and repelling. In particular, $Y\circ X^{-1}$ is uniformly expanding near its Julia set with respect to a complete orbifold length metric. 

Let $V \subset \MMM_P\union\{0\}$ be a simply-connected domain containing $0$ and $m_\circledast$.  The expanding property implies that for $r\in \N$ large enough, there exist domains $U_0, U_{m_\circledast}$ compactly contained in $V$, and univalent inverse branches of $R^r$ giving analytic isomorphisms $R_0^{-r}: V \to U_0$ and $R_{m_\circledast}^{-r}: V \to U_{m_\circledast}$ that fix $0$ and $m_\circledast$, respectively. 
The maps $R^{-r}_0, R^{-r}_{m_\circledast}$ are branches of the $r$th iterate of $X \circ Y^{-1}$; together with their domains, we have a so-called iterated function system.   By contraction, the fixed-points of these branches are unique.

In $\MMM_P$, the loci defined by the conditions $\ell(\tau) < \epsilon$ and $d(\pi_P(\tau), m_\circledast) < \delta$ contain neighborhoods of the end $0$ and the basepoint $m_\circledast$ of $\MMM_P$.   Define a sequence of points $\mu_n \in \MMM_P$ comprising an orbit under the pullback correspondence as follows. Set $\mu_0:=m_\circledast$.  Apply the branch  $R_0^{-r}: V \to U_0$ enough times so that the length of the systole becomes less than $\epsilon_1$. Now apply the branch $R_{m_\circledast}^{-r}$ enough times so that the resulting point lies in the neighborhood $d(\pi_P(\tau), m_\circledast) < \delta$. Now again apply the branch  $R_0^{-r}: V \to U_0$ enough times so that the length of the systole becomes less than $\epsilon_2$, etc. In this fashion, we obtain a finite orbit segment of the pullback correspondence in $\MMM_P$.  By conclusion $(I)$, the virtual endomorphism $\phi_f$ on $G_P$ is surjective.  Hence Proposition \ref{prop:shadowing} applies. We conclude that the orbit in $\MMM_P$ can be lifted to an orbit in $\TTT_P$ having the desired properties. \eb
\qed

\noindent{\bf Example.} Consider Example 2 from Section~\ref{section:FGA}.  As shown in \cite[\S~6.2]{bartholdi:nekrashevych:twisted}, the dendrite Julia set polynomial $f(z)=z^2+i$ satisfies the hypothesis of the previous theorem. Therefore all conclusions of the theorem follow, moreover, the following is true. Recall that $Y\circ X^{-1}$ is conjugate to the rational map $g(w) = (1-2/w)^2$, the Julia set of which is the entire Riemann sphere. Therefore, the forward orbit of any point in $\MMM_P$ under the pullback correspondence, which is the same as the backward orbit of $Y\circ X^{-1}$, is dense in $\MMM_P$. In other words, for any point $\tau \in \TTT_P$, the projection $\pi_P( \bigcup_{g \in G_P, n \in \N} \{\sigma_f^n(g\cdot \tau)\})$ of the union of $\sigma_f$-orbits of all points in $\TTT_P$, that are in the same fiber as $\tau$ over $\MMM_P$, is dense in $\MMM_P$. This situation is somewhat surprising as we know that every particular orbit $\sigma_f^n(g\cdot \tau)$ converges geometrically fast to the unique attracting fixed point $\tau_\circledast$, which corresponds to $f$. For this example, the statement of the conclusion (VI) in the last theorem can be further strengthened. We can find a point $\tau_0$ in the fiber of $m_\circledast$ that follows an arbitrary itinerary in the moduli space with an arbitrary precision before converging to the fixed point.

\bibliographystyle{math}
\bibliography{refs}

\def\cprime{$'$}
\begin{thebibliography}{BEKP}

\bibitem[BN1]{bartholdi:nekrashevych:twisted}
Laurent Bartholdi and Volodymyr Nekrashevych.
\newblock {Thurston equivalence of topological polynomials}.
\newblock {\em Acta Math.} {\bf 197}(2006), 1--51.

\bibitem[BN2]{bartholdi:nekrashevych:imgquad1}
Laurent Bartholdi and Volodymyr Nekrashevych.
\newblock {Iterated Monodromy Groups of Quadratic Polynomials, {I}}.
\newblock \url{http://arxiv.org/abs/math/0611177}, 2008.

\bibitem[BM]{bell:margalit:cohopfian}
Robert~W. Bell and Dan Margalit.
\newblock {Injections of {A}rtin groups}.
\newblock {\em Comment. Math. Helv.} {\bf 82}(2007), 725--751.

\bibitem[BLM]{blm:duke:1983}
Joan~S. Birman, Alex Lubotzky, and John McCarthy.
\newblock {Abelian and solvable subgroups of the mapping class groups}.
\newblock {\em Duke Math. J.} {\bf 50}(1983), 1107--1120.

\bibitem[BEKP]{bekp}
Xavier Buff, Adam Epstein, Sarah Koch, and Kevin Pilgrim.
\newblock {On {T}hurston's pullback map}.
\newblock In {\em Complex dynamics}, pages 561--583. A K Peters, Wellesley, MA,
  2009.

\bibitem[CFP]{cfp:fsr}
J.~W. Cannon, W.~J. Floyd, and W.~R. Parry.
\newblock {Finite subdivision rules}.
\newblock {\em Conform. Geom. Dyn.} {\bf 5}(2001), 153--196 (electronic).

\bibitem[CFPP]{cfpp:fsr-cve}
J.~W. Cannon, W.~J. Floyd, W.~R. Parry, and K.~M. Pilgrim.
\newblock {Subdivision rules and virtual endomorphisms}.
\newblock {\em Geom. Dedicata} {\bf 141}(2009), 181--195.

\bibitem[DH]{DH1}
A.~Douady and John Hubbard.
\newblock {A Proof of {T}hurston's Topological Characterization of Rational
  Functions}.
\newblock {\em Acta. Math.} {\bf 171}(1993), 263--297.

\bibitem[HP1]{kmp:ph:cxcii}
Peter Ha{\"{\i}}ssinsky and Kevin~M. Pilgrim.
\newblock {Finite type coarse expanding conformal dynamics}.
\newblock {\em Groups Geom. Dyn.} {\bf 5}(2011), 603--661.

\bibitem[HP2]{kmp:ph:expanding}
Peter Ha\"{\i}ssinsky and Kevin~M. Pilgrim.
\newblock {An algebraic characterization of expanding Thurston maps}.
\newblock \url{arXiv:1204.3214}, 2012. Submitted.

\bibitem[HK]{hubbard:koch:dm}
John Hubbard and Sarah Koch.
\newblock {An analytic construction of the Deligne-Mumford compactification of
  the moduli space of curves}.
\newblock \url{http://www.math.harvard.edu/~kochs/DMC.pdf}, 2012.

\bibitem[Hub]{hubbard:teich_vol1}
John~Hamal Hubbard.
\newblock {\em Teichm\"uller theory and applications to geometry, topology, and
  dynamics. {V}ol. 1}.
\newblock Matrix Editions, Ithaca, NY, 2006.
\newblock Teichm\"uller theory, With contributions by Adrien Douady, William
  Dunbar, Roland Roeder, Sylvain Bonnot, David Brown, Allen Hatcher, Chris
  Hruska and Sudeb Mitra, With forewords by William Thurston and Clifford
  Earle.

\bibitem[Iva1]{ivanov:book:subgroups}
Nikolai~V. Ivanov.
\newblock {\em Subgroups of {T}eichm\"uller modular groups}, volume 115 of {\em
  Translations of Mathematical Monographs}.
\newblock American Mathematical Society, Providence, RI, 1992.
\newblock Translated from the Russian by E. J. F. Primrose and revised by the
  author.

\bibitem[Iva2]{ivanov:mcg:1998}
Nikolai~V. Ivanov.
\newblock {Mapping class groups}.
\newblock \url{http://www.math.msu.edu/~ivanov/m99.ps}, 1998.

\bibitem[Kam]{kameyama:equivalence}
Atsushi Kameyama.
\newblock {The {T}hurston equivalence for postcritically finite branched
  coverings}.
\newblock {\em Osaka J. Math.} {\bf 38}(2001), 565--610.

\bibitem[Kel]{kelsey:schemes}
Gregory~A. Kelsey.
\newblock {Mapping schemes realizable by obstructed topological polynomials}.
\newblock {\em Conform. Geom. Dyn.} {\bf 16}(2012), 44--80.

\bibitem[Koc1]{koch:thesis}
Sarah Koch.
\newblock {\em Teichm\"uller theory and endomorphisms of $P^n$}.
\newblock PhD thesis, University of Provence, 2007.

\bibitem[Koc2]{koch:criticallyfinite}
Sarah Koch.
\newblock {TeichmŸller theory and critically finite endomorphisms}.
\newblock \url{http://www.math.harvard.edu/~kochs/endo.pdf}, 2012.

\bibitem[KPS]{koch:pilgrim:selinger:corlim}
Sarah Koch, Kevin Pilgrim, and Nikita Selinger.
\newblock {Limit sets of mapping class semigroups}.
\newblock Manuscript in preparation, 2012.

\bibitem[Lod]{lodge:thesis}
Russell Lodge.
\newblock {\em The boundary values of Thurston's pullback map}.
\newblock PhD thesis, Indiana University, 2012.

\bibitem[Mas]{masur:duke:1976}
Howard Masur.
\newblock {Extension of the {W}eil-{P}etersson metric to the boundary of
  {T}eichmuller space}.
\newblock {\em Duke Math. J.} {\bf 43}(1976), 623--635.

\bibitem[Nek1]{nekrashevych:combinatorics}
Volodymyr Nekrashevych.
\newblock {Combinatorics of polynomial iterations}.
\newblock In {\em Complex dynamics}, pages 169--214. A K Peters, Wellesley, MA,
  2009.

\bibitem[Nek2]{nekrashevych:expanding}
Volodymyr~V. Nekrashevych.
\newblock {Combinatorial models of expanding dynamical systems}.
\newblock \url{http://arxiv.org/pdf/0810.4936.pdf}, 2009.

\bibitem[Pil1]{kmp:canonical}
Kevin~M. Pilgrim.
\newblock {Canonical {T}hurston obstructions}.
\newblock {\em Adv. Math.} {\bf 158}(2001), 154--168.

\bibitem[Pil2]{kmp:cds}
Kevin~M. Pilgrim.
\newblock {Combination, decomposition, and structure theory for postcritically
  finite branched coverings of the two-sphere to itself}.
\newblock to appear, Springer Lecture Notes in Math., 2003.

\bibitem[Pil3]{kmp:tw}
Kevin~M. Pilgrim.
\newblock {An algebraic formulation of Thurston's characterization of rational
  functions}.
\newblock to appear, Annales de la Faculte des Sciences de Toulouse.

\bibitem[Sel1]{selinger:boundary}
Nikita Selinger.
\newblock {Thurston's pullback map on the augmented {T}eichm\"uller space and
  applications}.
\newblock {\em Invent. Math.} {\bf 189}(2012), 111--142.

\bibitem[Sel2]{selinger:euclidean}
Nikita Selinger.
\newblock {Topological characterization of canonical Thurston obstructions}.
\newblock arxiv \url{http://arxiv.org/abs/1202.1556}, 2012.

\bibitem[Wol]{wolpert:handbook}
Scott~A. Wolpert.
\newblock {The {W}eil-{P}etersson metric geometry}.
\newblock In {\em Handbook of {T}eichm\"uller theory. {V}ol. {II}}, volume~13
  of {\em IRMA Lect. Math. Theor. Phys.}, pages 47--64. Eur. Math. Soc.,
  Z\"urich, 2009.

\end{thebibliography}

\end{document}